\documentclass[reqno,11pt]{article}
\usepackage{amsmath,amssymb}
\usepackage{geometry}
\usepackage{amsmath, amsfonts,amsthm,amssymb,amsbsy,upref,color,graphicx,amscd,enumerate}
\usepackage[active]{srcltx}
\usepackage[utf8]{inputenc}
\usepackage{algorithm}
\usepackage{stmaryrd}
\usepackage{array}
\usepackage{algorithmicx}
\usepackage{algpseudocode}
\usepackage{graphicx}
\usepackage{color}

\usepackage[table,xcdraw,svgnames]{xcolor}  
\usepackage[colorlinks]{hyperref}  
\hypersetup{
	allbordercolors=.,
	citecolor=NavyBlue,
	linkcolor=DarkRed,
	urlcolor=NavyBlue
}

\usepackage{multirow}
\graphicspath{{./metapost/}{./pics/}}

\usepackage{enumitem}

\numberwithin{equation}{section}
\theoremstyle{plain}
\newtheorem{thm}{Theorem}[section]
\newtheorem{lemma}{Lemma}
\newtheorem{cor}[thm]{Corollary}
\newtheorem{prop}[thm]{Proposition}
\newtheorem{deff}[thm]{Definition}
\theoremstyle{definition}
\newtheorem{rem}[thm]{Remark}
\newtheorem{conj}{Conjecture}

\renewcommand{\Bbb}{\mathbb}
\newcommand{\bo}[1]{{\bf #1}}

\newcommand{\Red}[1]{{\color{DarkRed}#1}}
\definecolor{awesome}{rgb}{1.0, 0.13, 0.32}
\definecolor{coolblack}{rgb}{0.0, 0.18, 0.39}
\definecolor{darkcerulean}{rgb}{0.03, 0.27, 0.49}
\newcommand{\Label}[1]{{\color{black}#1}}

\makeatletter
\DeclareFontFamily{U}{tipa}{}
\DeclareFontShape{U}{tipa}{m}{n}{<->tipa10}{}
\newcommand{\arc@char}{{\usefont{U}{tipa}{m}{n}\symbol{62}}}%

\newcommand{\arc}[1]{\mathpalette\arc@arc{#1}}

\newcommand{\arc@arc}[2]{%
  \sbox0{$\m@th#1#2$}%
  \vbox{
    \hbox{\resizebox{\wd0}{\height}{\arc@char}}
    \nointerlineskip
    \box0
  }%
}
\makeatother

\title{Volume computation for Meissner polyhedra and applications}

\author{Beniamin Bogosel}

\begin{document}
	
\maketitle

\begin{abstract}
The volume of a Meissner polyhedron is computed in terms of the lengths of its dual edges. This allows to reformulate the Meissner conjecture regarding constant width bodies with minimal volume as a series of explicit finite dimensional problems. A direct consequence is the minimality of the volume of Meissner tetrahedras among Meissner pyramids. 
\end{abstract}

\section{Introduction}

Certain convex shapes have the remarkable property that any two parallel tangent or supporting planes which contain the shape between them are at a fixed distance apart. The ball is an obvious example, but there are infinitely many more. In dimension two the famous Reuleaux triangle (the intersection of three unit disks centered at the vertices of a unit equilateral triangle) has constant width, together with the whole class of Reuleaux polygons. See \cite[Chapter 7]{yaglom-boltjanskii} for an introduction regarding shapes of constant width. For simplicity, all constant width shapes have unit diameter in the following.

The Reuleaux triangle is a particular shape of constant width, since it solves various optimization problems in the class of constant width:
\begin{itemize}[noitemsep,topsep=0pt]
	\item it minimizes the area \cite{Blaschke1915}, \cite{Lebesgue}, \cite{Chakerian}, \cite{yaglom-boltjanskii};
	\item it minimizes the inradius \cite[Chapter 7]{yaglom-boltjanskii}, or equivalently, it maximizes the circumradius;
	\item it maximizes the Cheeger constant \cite{henrot-lucardesi}, \cite{bogosel-cheeger}.
\end{itemize}
In dimension three fewer results are known regarding extremal shapes of constant width. The minimality of the inradius and the maximality of the circumradius are achieved by constant width shapes containing a unit tetrahedron \cite[Section 14.3]{bodies_of_constant_width}. In \cite{anciaux_guilfoyle} it is proved that any body of constant width minimizing the volume resembles the Meissner tetrahedron in the sense that any diameter has one endpoint corresponding to a singular part of the boundary. Procedures for producing shapes of constant width in arbitrary dimension were given in \cite{lachand-robert-oudet} and analytical parametrizations for three dimensional constant width bodies were proposed in \cite{bayen-LR-Oudet}.

Nevertheless, there is are three dimensional bodies which are conjectured to minimize the volume. These bodies are called \bo{Meissner tetrahedra} \cite{Meissner-Schilling} and are constructed as follows:
\begin{itemize}[noitemsep,topsep=0pt]
	\item Intersect four unit balls centered at the vertices of a regular tetrahedron of unit diameter. The body obtained is called a Reuleaux tetrahedron, but it does not have constant width (see \cite[Section 8.2]{bodies_of_constant_width} for more details)
	\item One side is chosen among each pair of opposite sides and a \emph{smoothing procedure} is performed. This is described in \cite[Chapter 7]{yaglom-boltjanskii} or \cite[Section 8.3]{bodies_of_constant_width}. The smoothing procedure can give rise to two types of Meissner tetrahedra: either all edges coming from a vertex are smoothed or all edges adjacent to one of the faces are smoothed.
\end{itemize}
It is conjectured by Bonnesen and Fenchel \cite{Bonnesen-Fenchel} that Meissner tetrahedra minimize the volume. More historical and bibliographical aspects related to this problem are given in \cite{kawohl-webe}. Numerical simulations presented in \cite{AntunesBogosel22} further suggest that the Meissner tetrahedra are indeed minimizers for the volume. The conjecture is still open today and we state it below:
\begin{conj}\label{conj:meissner}
	The Meissner tetrahedra minimize the volume among all three dimensional bodies with fixed constant width.
\end{conj}

Many proofs of the analogue two dimensional result rely on showing that any Reuleaux polygon which is not the Reuleaux triangle is not a minimizer for the area. The density of Reuleaux polygons in the class of two dimensional shapes of constant width implies the result. 

Until recently, a similar class of discrete constant width shape was missing in dimension three. In \cite{sallee_polytopes} the class of Reuleaux polytopes was introduced, defined using intersection of balls. Although these polytopes can approximate arbitrarily well all bodies of constant width, they have a drawback: they do not have constant width themselves. This was remedied recently in \cite{montejano} where the authors show how to obtain constant width bodies in dimension three starting from a Reuleaux polytope. The resulting bodies are called Meissner polyhedra, since they generalize the Meissner tetrahedron. More details regarding these bodies, including their connections with extremal sets of diameter one and their density in the class of shapes of constant width are presented in \cite{meissner_hynd}. It turns out that Meissner polyhedra give a natural context for studying Conjecture \ref{conj:meissner}, since the Meissner tetrahedron naturally belongs to all classes of Meissner polyhedra which have an upper bound on the number of vertices. 

The study of Meissner polyhedra gives rise to complex combinatorial aspects and connections with graph theory. The graphs behind the Reuleaux and Meissner polyhedra are studied in \cite{meissner_graphs}. In particular, on associated the Github page 
\begin{center}
	\href{https://github.com/mraggi/ReuleauxPolyhedra}{\nolinkurl{https://github.com/mraggi/ReuleauxPolyhedra}}
\end{center}
the graphs giving rise to Meissner polyhedra up to $14$ vertices are presented. These graphs are used for creating some of the illustrations in this paper. Recently, in \cite{hynd-vol-per} the surface area and volume of Reuleaux and Meissner polyhedra are computed, giving new tools for attacking Conjecture \ref{conj:meissner}.

The purpose of this paper is to present observations and computations which lead to a simple formula for the surface area and volume of Meissner polyhedra. Section \ref{sec:meissner} presents details regarding the construction of Meissner polyhedra. Section \ref{sec:area-computations} presents the detailed computation of the area and volume of Meissner polyhedra. In particular, Conjecture \ref{conj:meissner} is reduced to a maximization problem involving a completely explicit two dimensional function and the lengths of pairs of dual edges in a Meissner polyhedron. In Section \ref{sec:pyramids} the particular case of Meissner pyramids, generalizing the Meissner tetrahedron is discussed, establishing that the tetrahedron is the body of constant width with minimal volume among all pyramids of constant width. In Section \ref{sec:existence} we show that Conjecture \ref{conj:meissner} may be reduced to a series of finite dimensional problems for which the Meissner tetrahedra are natural solution candidates.

Illustrations shown in the paper were produced using Metapost or Matlab. The codes associated to \cite{spherical-design} were used to create drawings from spherical geometry.

\section{Meissner Polyhedra}
\label{sec:meissner}

All constant width bodies considered in the following have unit diameter. The Meissner Polyhedra are introduced in \cite{montejano} and are the 3D analogue of Reuleaux polygons. In \cite{montejano} it is mentioned that these polyhedra are dense in the class of 3D constant width bodies using arguments based on \cite{sallee_polytopes}. The density of Meissner polyhedra in the class of three dimensional constant width bodies was revisited in \cite{meissner_hynd} where a detailed description is given regarding the construction of these bodies. Moreover, a detailed description of the computation of the volume and surface area for the Meissner tetrahedra is given. The surface area and the volume of Meissner polyhedra is computed in \cite{hynd-vol-per}. In this paper an alternate computation is proposed, resulting in a simple formula depending only on the lengths of the pairs of dual edges.

 The volume of a  constant width body $K$ is related to the surface area using the Blaschke Formula
\begin{equation}\label{eq:blaschke}
 |K| = \frac{1}{2} |\partial K| -\frac{\pi}{3}. 
 \end{equation}
Therefore, to study the bodies having minimal volume it is enough to find constant width bodies having minimal surface areas. Moreover, since the Meissner polyhedra are dense in the class of constant width bodies if one proves that an arbitrary Meissner polyhedron has surface area bigger than the Meissner tetrahedron, the Meissner conjecture, i.e. the Blaschke-Lebesgue problem in dimension three is proved.

 In the following, we use the notation $B(X)$ to denote the intersection of all balls having unit radius with centers in $X$. Let $X\subset \Bbb{R}^3$ be a finite set of $m \geq 4$ points having diameter $1$ and suppose that $X$ is extremal, i.e. the number of diametric pairs of $X$ is maximal. In \cite[Theorem 3.4]{meissner_hynd} is recalled that a set is extremal if and only if there are $2m-2$ diametric pairs among points of $X$. Furthermore $X$ is the set of vertices of $B(X)$, i.e. the points where at least three spherical faces meet. In a similar manner, an edge of $B(X)$ is the intersection of two adjacent spherical faces of $\partial B(X)$. Following \cite{sallee_polytopes}, \cite{disk-polygons} and \cite{meissner_hynd} we can define Reuleaux polyhedra.
 
 \begin{deff}\bo{(Reuleaux polyhedra)}\label{eq:def-reuleaux-poly}
 	Given $X \subset \Bbb{R}^3$ a finite set of $m\geq 4$ points which is an extremal finite set of unit diameter, the associated Reuleaux polyhedron is $R = B(X)$.
 	
 	The resulting ball polytope does not have constant width. Nevertheless, any body of constant width can be approximated arbitrarily well using Reuleaux polyhedra. 
 \end{deff}

  Let $X$ be an extremal set of diameter one. Reuleaux polyhedra have natural definitions for vertices, edges and faces, which are analogue to classical polyhedra. In particular, vertices are points $x$ belonging to at least two diameters. Faces are spherical portions of $\partial B(X)$, each face being opposite to a vertex which is the center of the respective sphere. Two adjacent faces meet along an edge, which is a circular arc. See \cite{meissner_hynd} for a detailed discussion.
  
\begin{figure}
  	\centering 
  	\includegraphics[width=0.4\textwidth]{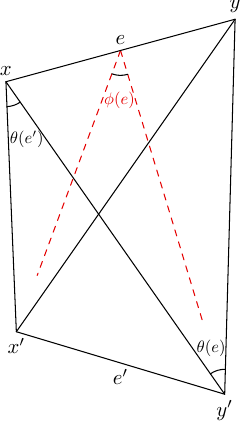}
  	\caption{Configuration of a pair of dual edges $(e,e')$. All other four edges in the corresponding tetrahedron have unit length. $\theta(e)$ denotes the spherical length of the edge $e$ and $\phi(e)$ represents the dihedral angle corresponding to edge $e$ in the tetrahedron determined by $(e,e')$.}
  	\label{fig:dual-edges}
  \end{figure}

   If $x,y$ are the endpoints of an edge $e$ of $B(X)$ then there is a unique edge $e'$ of $B(X)$ with endpoints $x',y'$ such that $|x-x'|=|x-y'|=|y-x|=|y-y'|=1$. See Figure \ref{fig:dual-edges} for an illustration. The pair $(e,e')$ is a pair of \emph{dual edges} of $B(X)$. See \cite[Section 3]{meissner_hynd} for more details.
  
  The Meissner polyhedra are defined as follows. Consider an extremal diameter one set $X\subset \Bbb{R}^3$ having $m\geq 4$, which according to \cite{meissner_hynd} has $m-1$ pairs of dual edges
  \begin{equation}\label{eq:dual-edges} (e_1,e'_1),...,(e_{m-1},e'_{m-1}).
  \end{equation}
  Meissner bodies were introduced in \cite{montejano} and we give the following equivalent definition according to \cite{meissner_hynd}.
  
\begin{figure}
	\centering
	\includegraphics[height=0.45\textwidth]{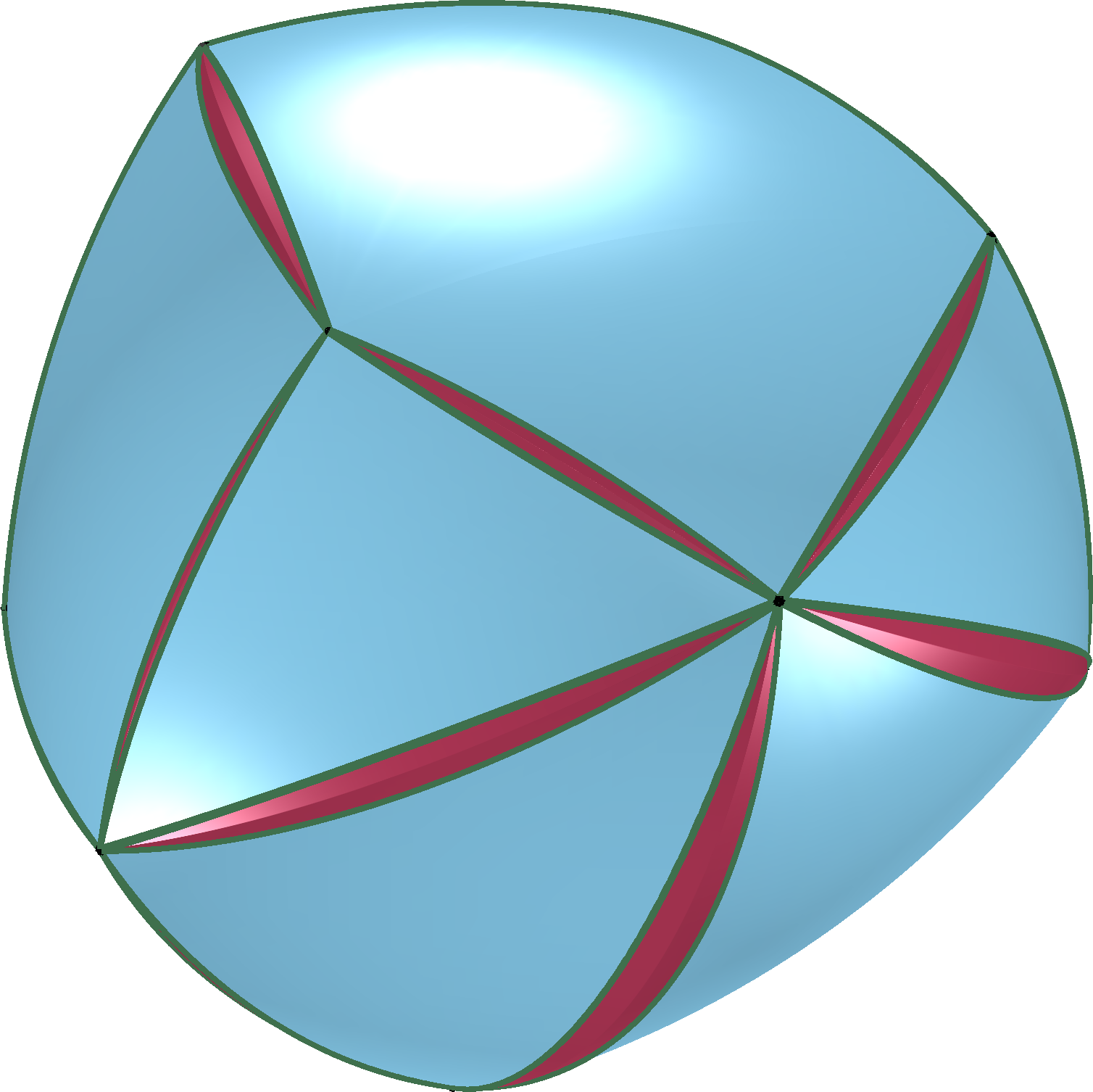}\quad
	\includegraphics[height=0.45\textwidth]{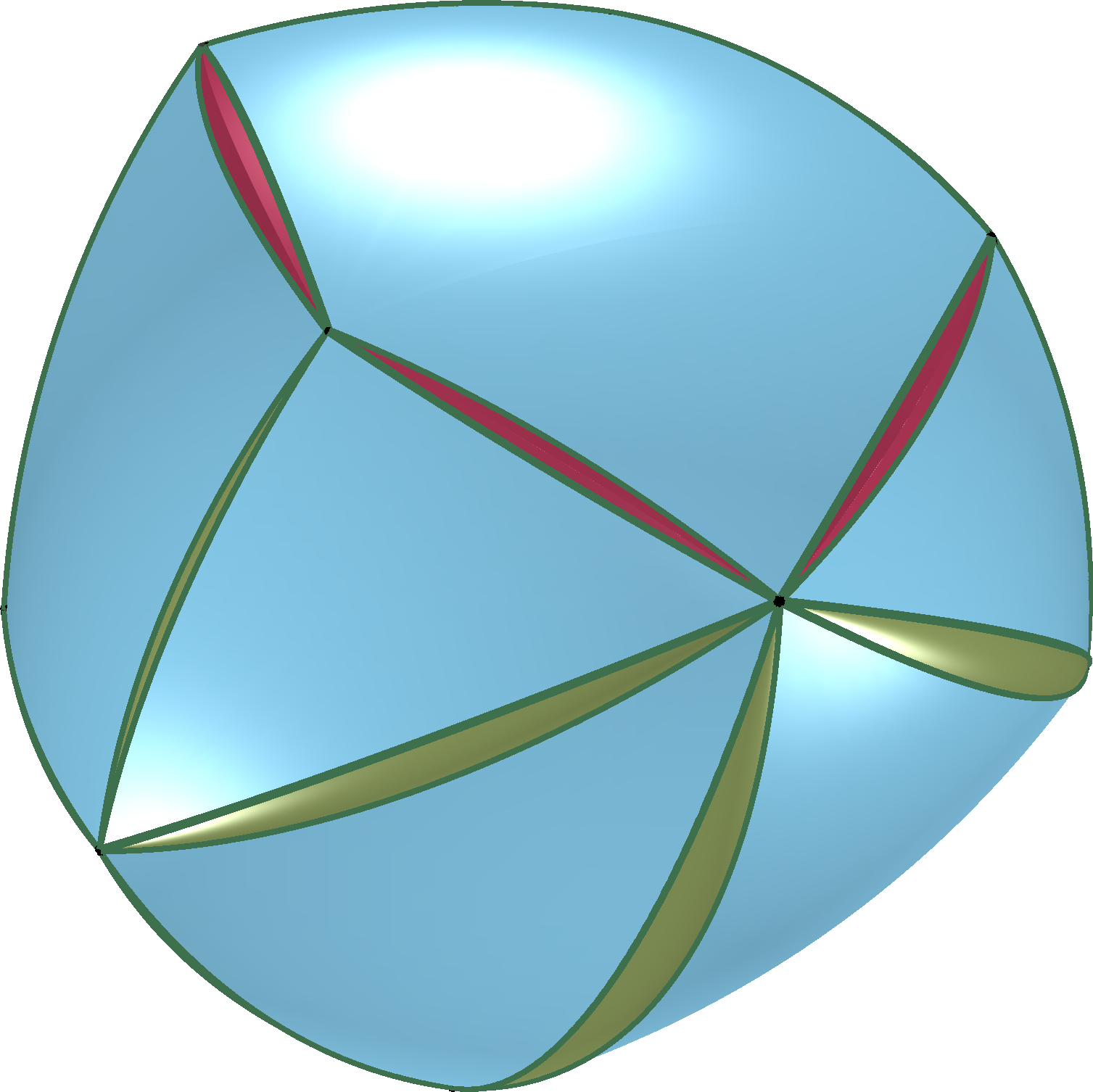}
	\caption{An example of Reuleaux polyhedron with $10$ vertices (left) and the associated Meissner polyhedron (right) obtained by smoothing some of the edges. Non-smoothed edges are represented in red, smoothed edges are represented with green and spherical geodesic polygons contained in the faces are colored in blue.}
	\label{fig:reuleaux-vs-meissner} 
\end{figure}

  \begin{deff}\bo{(Meissner polyhedra)}\label{def:meissner}
    Consider $X \subset \Bbb{R}^3$ a finite set of $m\geq 4$ points which is an extremal finite set of unit diameter having the pair of dual edges given by \eqref{eq:dual-edges}. The convex body
	\begin{equation}\label{eq:def_meissner}
	B(X\cup e_1\cup ... \cup e_{m-1})
	\end{equation}
	is a Meissner polyhedron based on $X$. Every Meissner polyhedron is a body of constant width.
  \end{deff}
    Note that it consists on intersection of balls with centers in $X$ and on one edge among every pair of dual edges. For every extremal finite set of unit diameter having $m$ points there exist $2^{m-1}$ choices of Meissner polyhedra which can be constructed. Among these choices there is one having minimal volume. This choice will be made more precise in the next section. The fact that this construction produces bodies of constant width is well established and proved in detail in \cite{montejano} and \cite{meissner_hynd}. An example of Reuleaux polyhedron and the associated Meissner polyhedron is shown in Figure \ref{fig:reuleaux-vs-meissner}. This example is taken from the database associated to the paper \cite{meissner_graphs}.
  
  Consider the following elements associated to an extremal diameter $1$ set:
  \begin{itemize}[noitemsep,topsep=0pt]
  	\item A generic pair of dual edges from \eqref{eq:dual-edges} will be denoted by $(e,e')$.
  	\item For each edge $e$ consider the associated angle $\theta(e)$ made by $e$ at the center of a unit ball which determines $e$.  More precisely, the endpoints of $x,y$ of $e$ are put on a disk of radius $1$ and the angle at the center is measured. In particular $\theta(e)$ measures the geodesic distance from the endpoints of $e$ on a unit sphere. Since the vertices form a set of diameter $1$ we must have $|x-y|\leq 1$ which implies $\theta(e) \in [0,\pi/3]$. {Moreover, if $\theta(e)=\theta(e')=\pi/3$ then $x,y,x',y'$ must coincide with the vertices of the regular tetrahedron.} We have the explicit formula $\sin (\theta(e)/2)= |x-y|/2$.
  	\item For each edge $e$ consider the tetrahedron determined by the vertices in the dual pair $(e,e')$ shown in Figure \ref{fig:dual-edges}. In this tetrahedron, edges $e,e'$ are orthogonal since all remaining edges have unit length. We denote by $\phi(e)$ the dihedral angle of this tetrahedron associated to the edge determined by $e$. A simple computation involving elementary trigonometry implies the following relation between $\theta(e),\theta(e')$ and $\phi(e)$:
  	\[ \sin \frac{\phi(e)}{2} = \frac{\sin \frac{ \theta(e')}{2}}{\cos \frac{\theta(e)}{2}}.\]
  	The elements of this tetrahedron are completely determined by $\theta(e), \theta(e')$.
  	\item For each vertex $x$ consider the opposite face $\tau(x)$ determined by vertices $x_1,...,x_k$, $k \geq 3$ situated at unit Euclidean distance from $x$. The extremal diameter one sets generating Meissner polyhedra have the self-dual property: $x \in \tau(y) \Longleftrightarrow y \in \tau(x)$. This notation for opposite faces is used for Meissner or Reuleaux polyhedra.
  \end{itemize}

At each pair of dual edges $(e,e')$ in a Meissner polyhedron two types of surfaces appear. If points of $e$ are chosen as centers of balls in the intersection determining the Meissner polyhedron $M$ then the surface of $M$ near $e'$ consists of a spindle surface, i.e. a surface obtained by rotating an arc of circle of radius $1$ around a symmetry axis. We say that such an edge is \emph{smoothed}. More details regarding the geometry of this configurations are given in \cite[Section 4.1]{meissner_hynd} and in the following.

\begin{figure}
	\centering
	\includegraphics[width=0.3\textwidth]{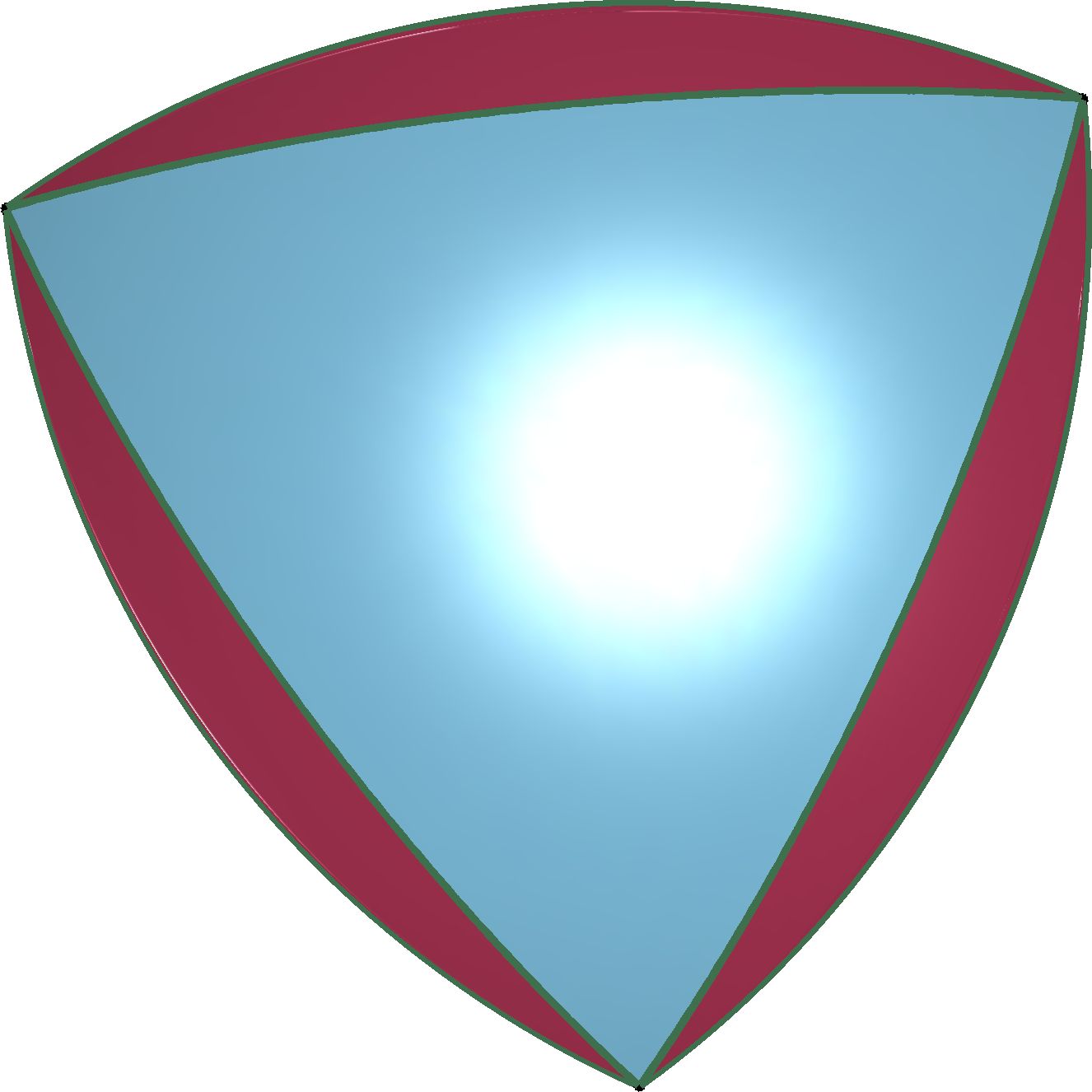}\quad
	\includegraphics[width=0.3\textwidth]{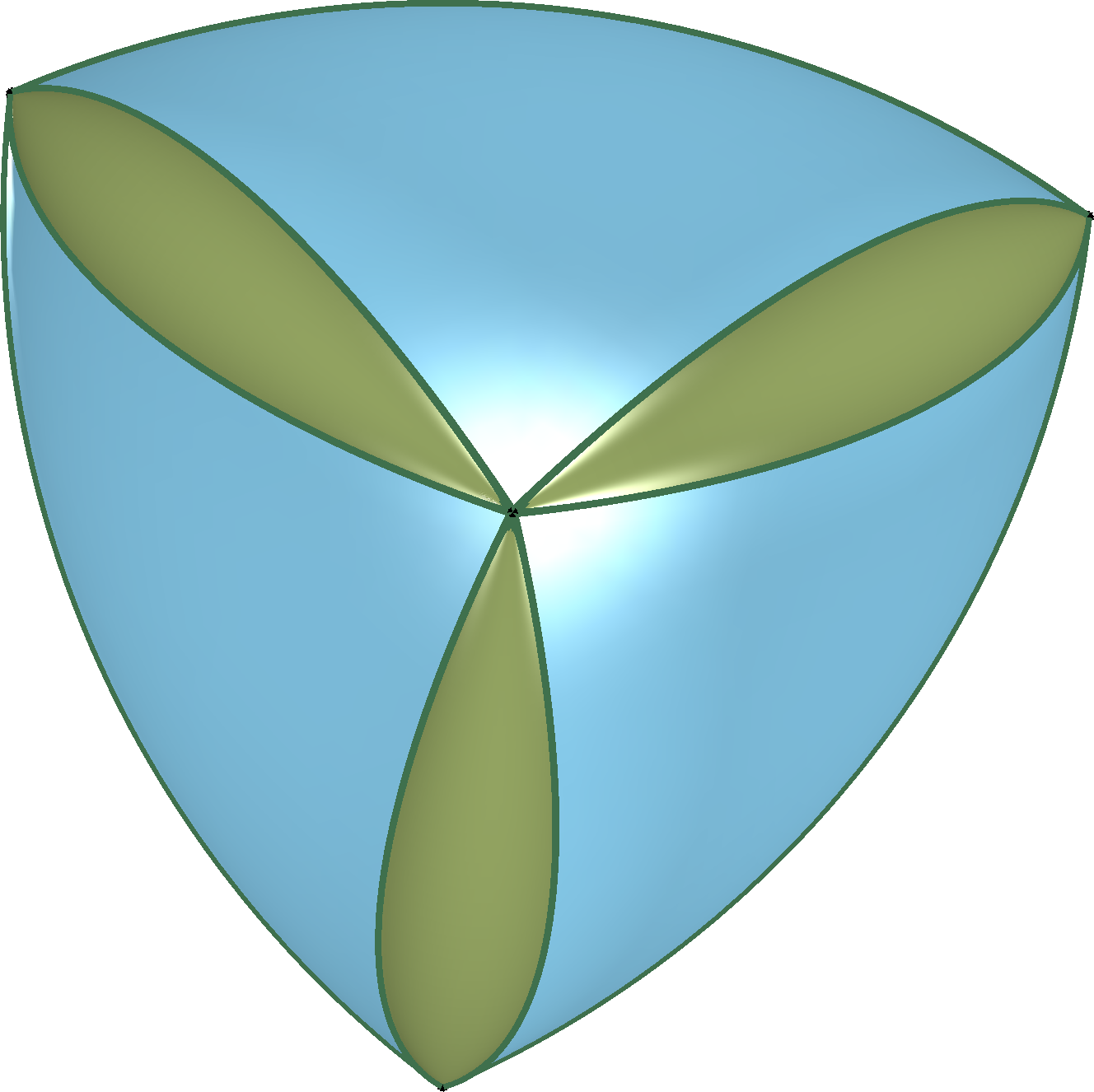}\quad
	\includegraphics[width=0.3\textwidth]{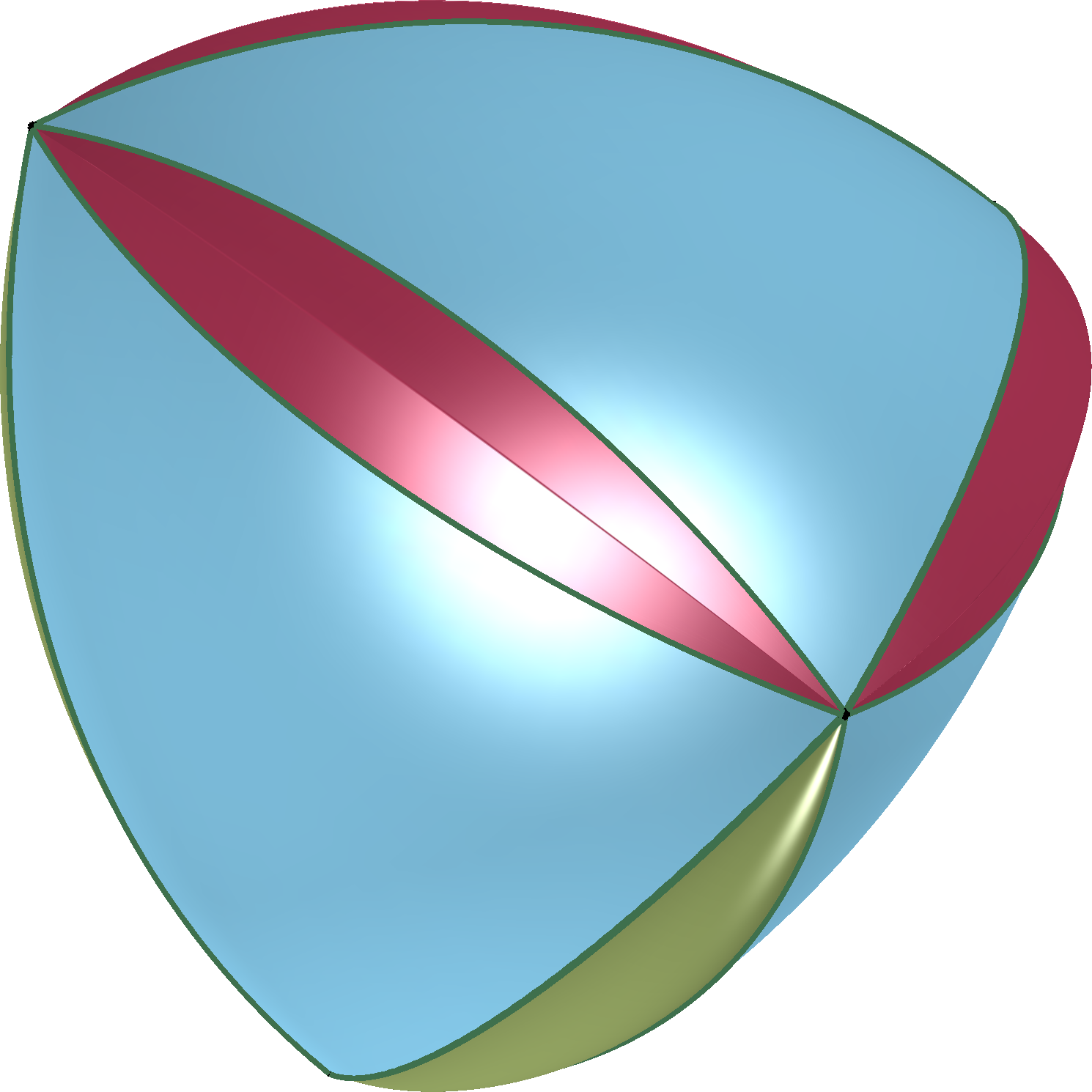}

	\medskip
	
	\includegraphics[width=0.3\textwidth]{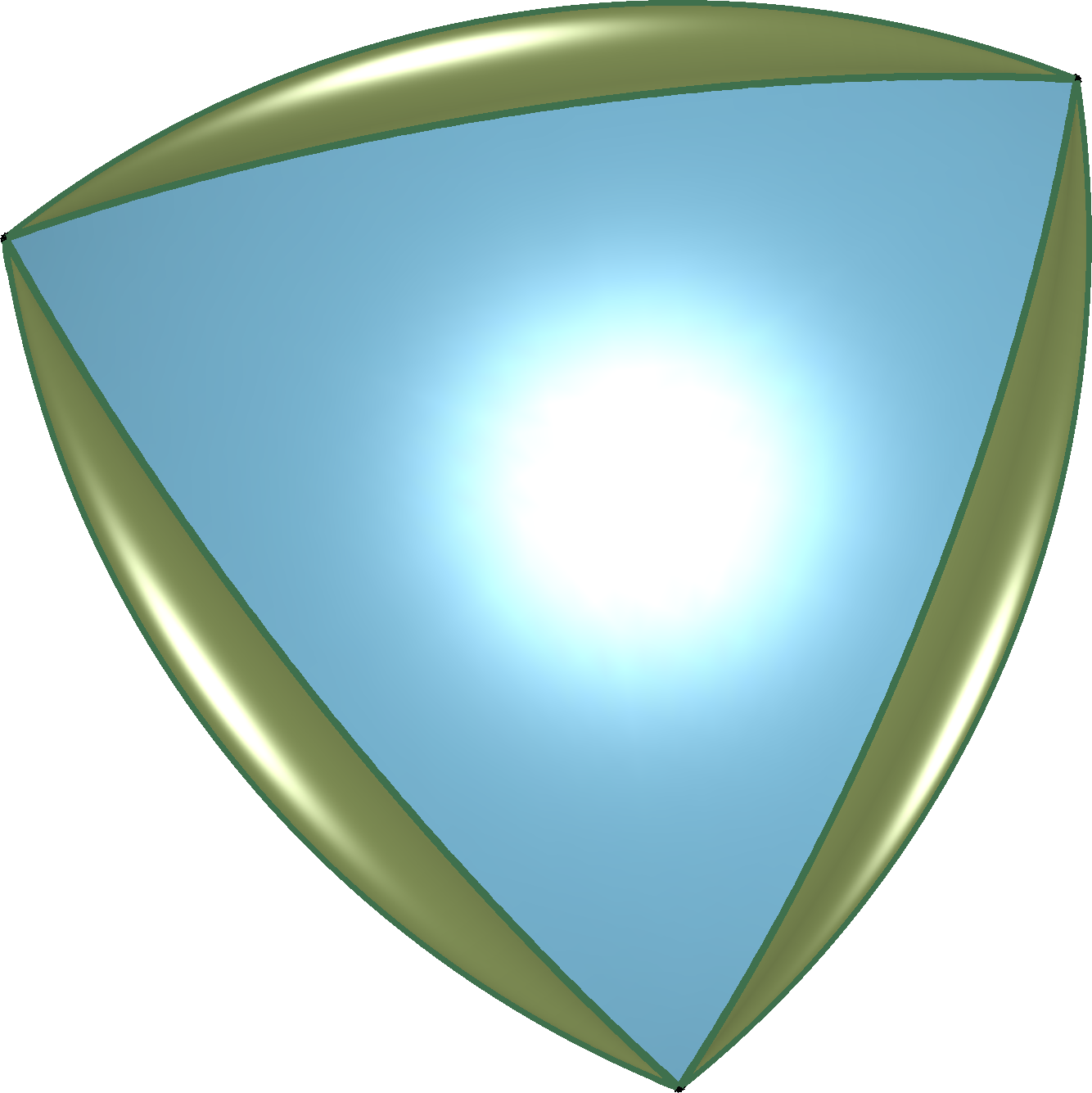}\quad
	\includegraphics[width=0.3\textwidth]{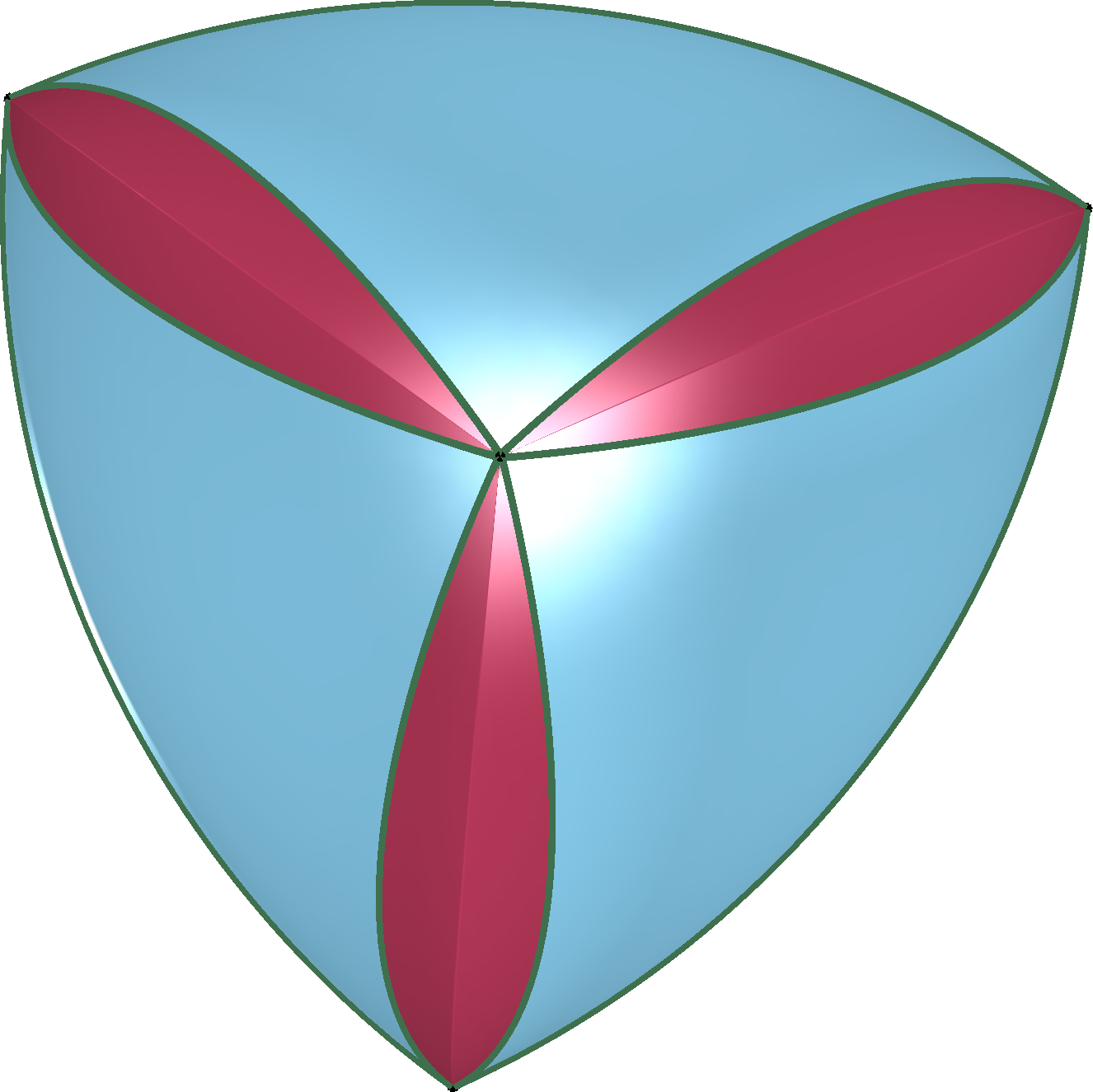}\quad
	\includegraphics[width=0.3\textwidth]{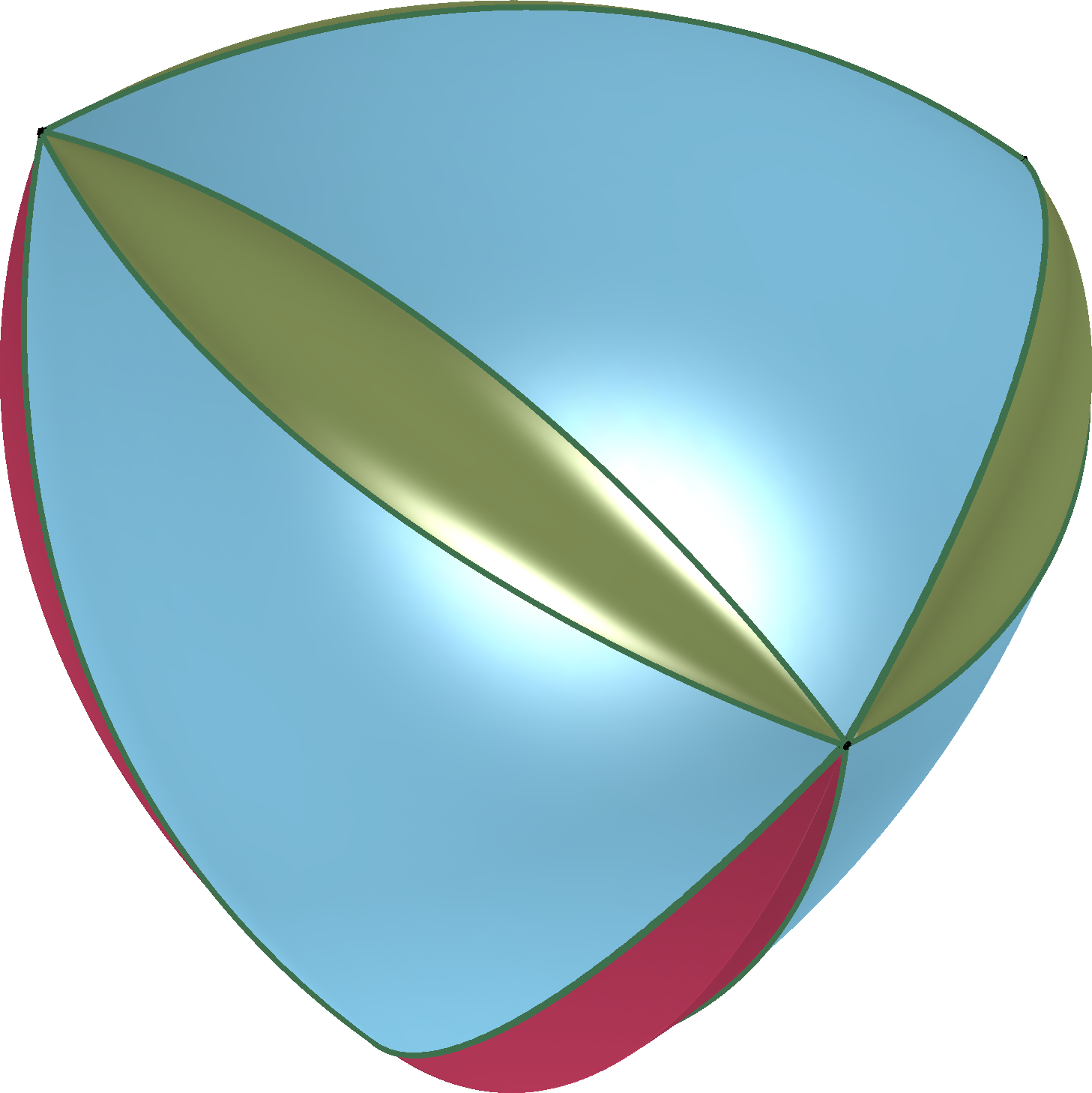}
	\caption{Meissner tetrahedra of the two types. Either three edges having a common vertex are smoothed (first row) or three edges adjacent to a common face are smoothed (second row). Spherical parts are represented in blue, wedge surfaces are shown in red and spindle surfaces in green.}
	\label{fig:tetraM}
\end{figure}
  
\begin{figure}
	\centering 
	\includegraphics[width=0.4\textwidth]{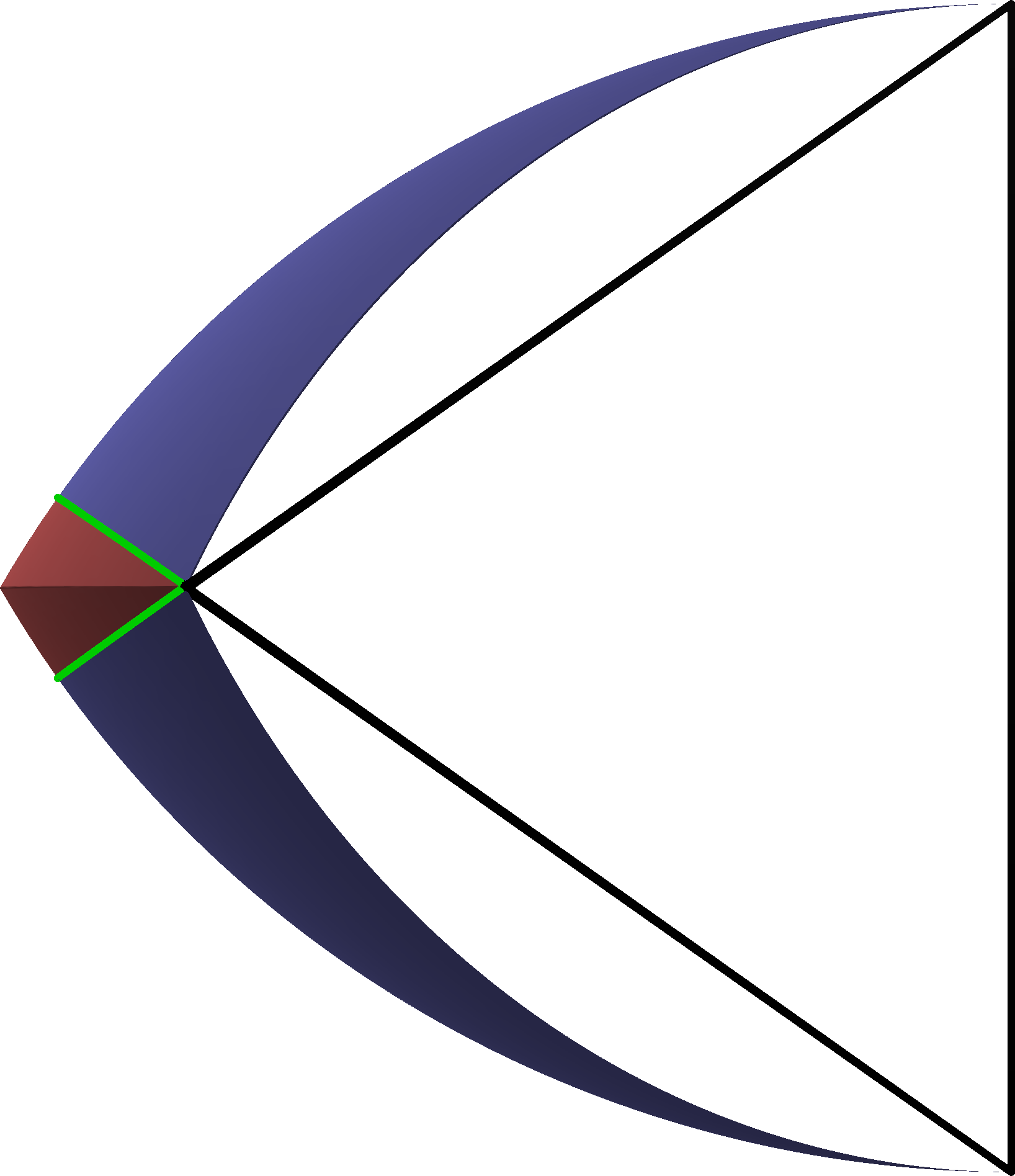}\quad
	\includegraphics[width=0.4\textwidth]{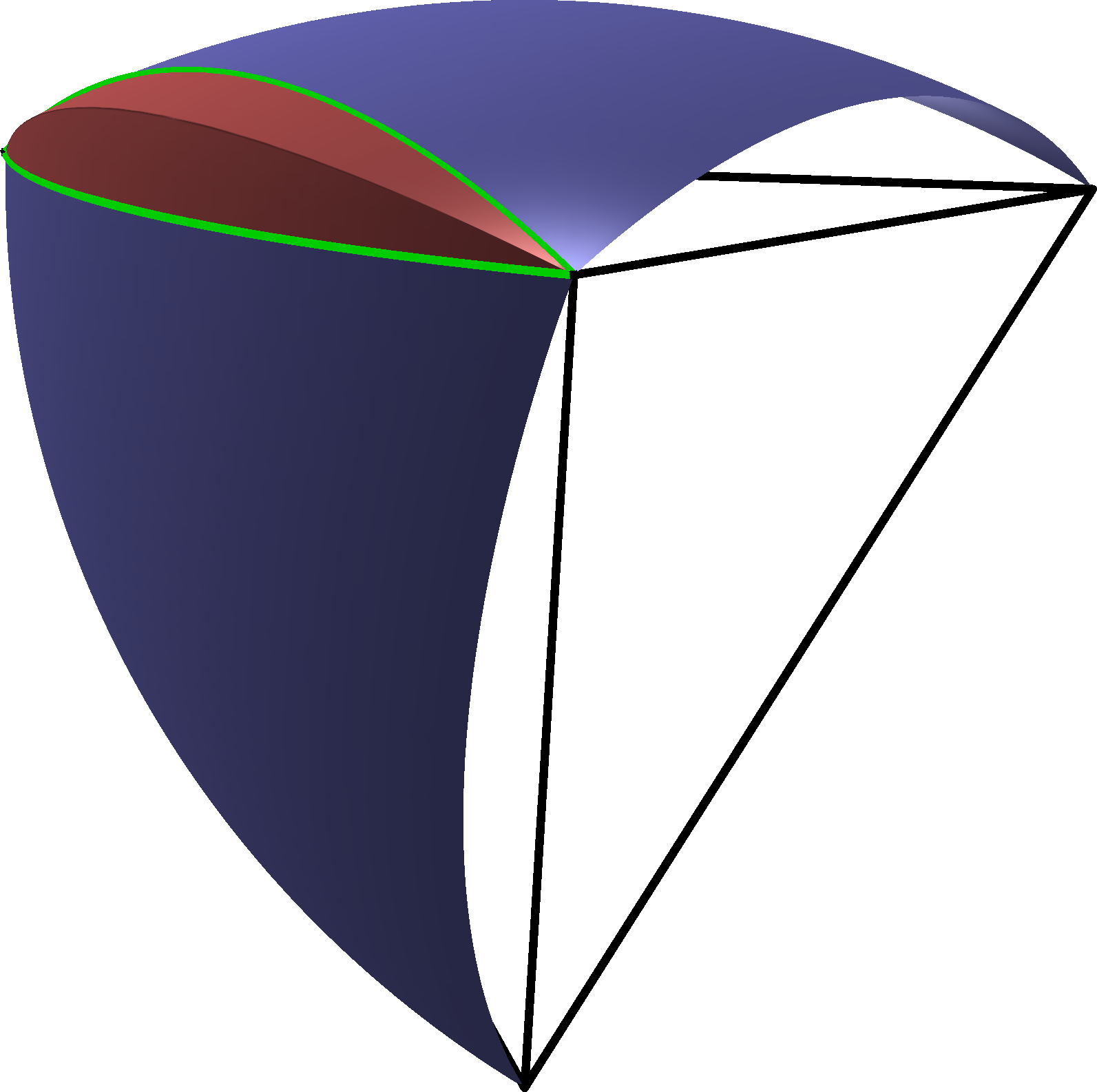}\quad
	\caption{Wedge surface associated to an edge $e$ with vertices $x,y$ of a Meissner polyhedron. The two adjacent spherical faces centered at $x',y'$ endpoints of the dual edge $e'$ are continued until they intersect across the geodesics connecting $x,y$ on these faces. }
	\label{fig:wedges}
\end{figure}

\begin{deff}\label{def:wedge}
{\bf Wedge surface.} Consider $(e,e')$ a pair of dual edges in the Meissner polyhedron $M$. Endpoints $x,y$ of the edge $e$ belong to the intersection of two spheres of radii $1$ centered at $x',y'$, the endpoints of $e'$. Denote by $\arc{xy}$ the small arc associated to $e$ in the circle $\partial B(x')\cap \partial B(y')$. Consider $g_{x'},g_{y'}$ the geodesic arcs joining $x,y$ on $\partial B(x'),\partial B(y')$, respectively. The wedge $W(e)$ at edge $e$ is defined as the region bounded by $\arc{xy}$ and $g_{x'}$ in $\partial B(x')$ and $ \arc{xy}$ and $g_{y'}$ in $\partial B(y')$. It is obvious that $W(e)$ is symmetric about $e$. See Figure \ref{fig:wedges} for an illustration.
\end{deff}

\begin{figure}
	\centering 
	\includegraphics[height=0.4\textwidth]{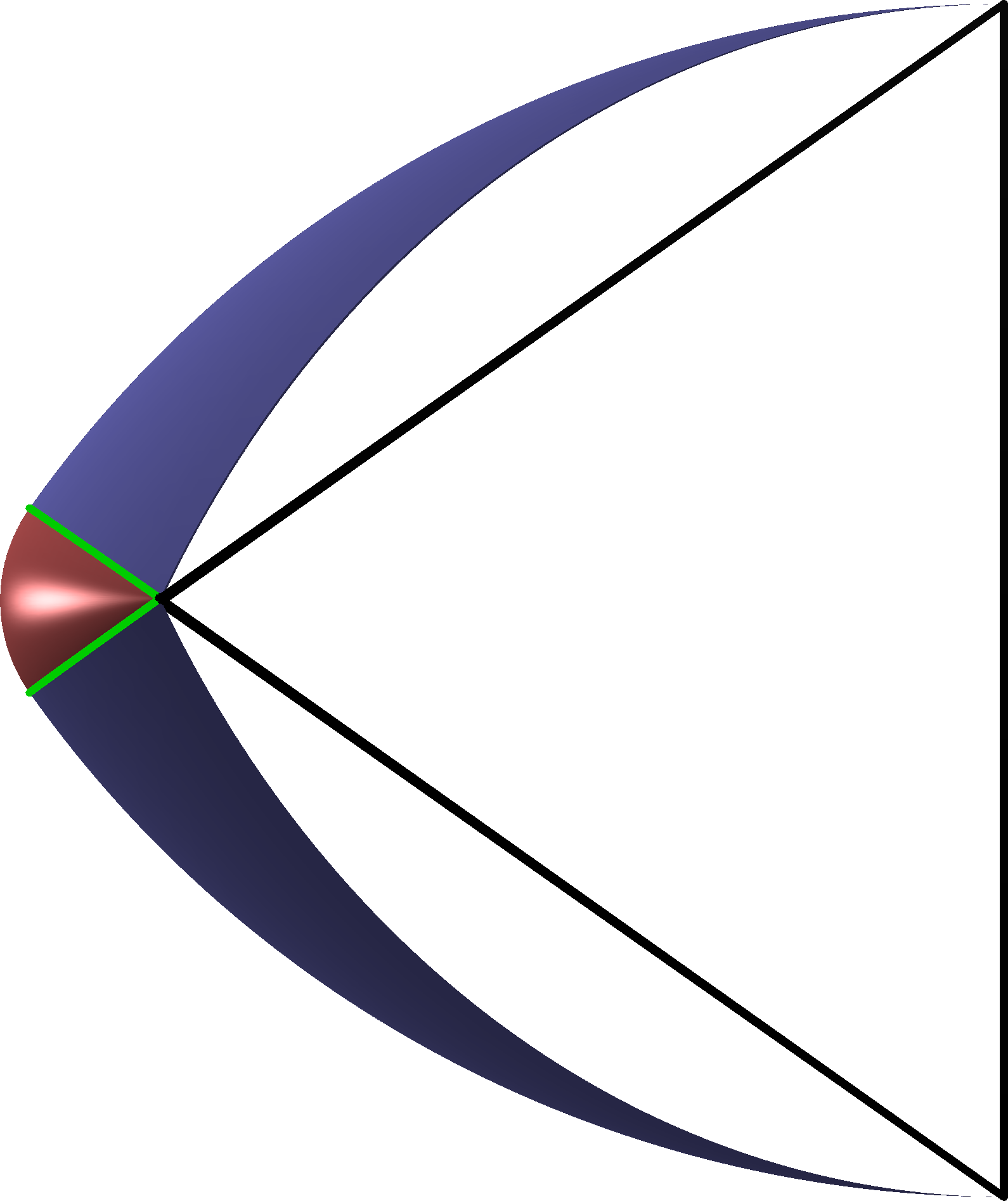}\quad 
	\includegraphics[height=0.4\textwidth]{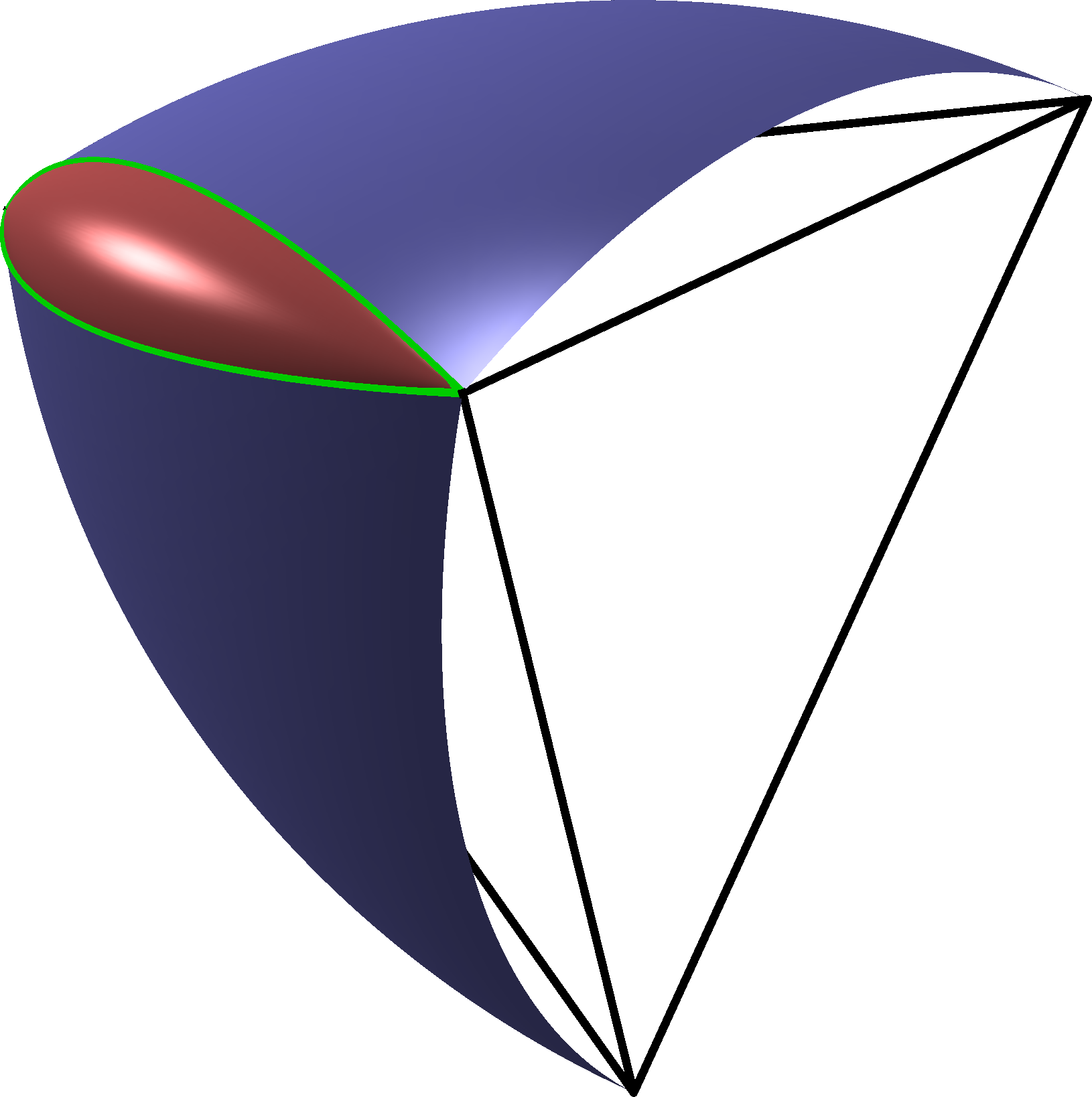}
	\caption{Spindle surface associated to an edge $e$ with vertices $x,y$ of a Meissner polyhedron. The geodesics connecting $x,y$ on the adjacent spherical faces are rotated around $xy$. }
	\label{fig:spindle}
\end{figure}

\begin{deff}\label{def:spindle}
	{\bf Spindle surface.} Consider $(e,e')$ a pair of dual edges in $M$. Using the notation defined previously, the geodesic arcs $g_{x'},g_{y'}$ are constructed around the edge $e$ bounded by vertices $x$ and $y$ on spheres $\partial B(x'),\partial B(y')$. Both $g_{x'}$ and $g_{y'}$ are geodesic arcs in sphere of radius $1$, they are arcs of circles of radius $1$ having length $\theta(e)$. The spindle surface $S(e)$ is the surface obtained by rotating $g_{x'}$ towards $g_{y'}$ around the axis $xy$. It is part of a surface of revolution determined by the dihedral angle $\phi(e)$ associated to edge $xy$ in the tetrahedron $x,y,x',y'$. See Figure \ref{fig:spindle} for an illustration. 
\end{deff}

\begin{rem}
	In \cite{montejano} spindle surfaces are called wedges, however this does not reflect the meaning of the word \emph{wedge} which is a triangular shaped tool, implying the existence of an angle. In \cite{meissner_hynd} the wedge surfaces are called \emph{silver surfaces}.
\end{rem}

From the previous definitions of wedge and spindle surfaces, it is apparent that their surface areas depend only on the parameters $\theta(e),\theta(e')$ characterizing the tetrahedron $x,y,x',y'$. Moreover, if $\theta(e)=\theta(e')$, i.e. the tetrahedron is symmetric, choosing $(e',e)$ instead of $(e,e')$ in \eqref{eq:def_meissner} will not change the area or volume of the Meissner polyhedron $M$.

If $m=4$ then the only extremal set of diameter $1$ consists of the vertices of a regular tetrahedron of edge length $1$. In this case, the Meissner polyhedra coincide with one of the two Meissner tetrahedra \cite[Chapter 7]{yaglom-boltjanskii}, \cite[Section 8.3]{bodies_of_constant_width}. Note that the $2^3$ choices of smoothing one edge among the three pair of opposite edges only produces two types of Meissner tetrahedra. Either all edges starting from a common vertex are smoothed or the edges adjacent to a given face are smoothed. See Figure \ref{fig:tetraM} for an illustration. 

The two Meissner tetrahedra have the same volume and the same area. As already underlined in \cite[Section 4]{meissner_hynd}, for a given extremal finite set of diameter one, the corresponding Meissner polyhedra differ only in the exterior dihedral angles associated to edges $e=xy,e'=x'y'$ in the tetrahedron $x,y,x',y'$. Indeed, all other regions of $\partial M$ are portions of spheres centered at the vertices of $M$. Moreover, wedge or spindle surfaces in $\partial M$ intersect only at vertices. Since the vertices of a Meissner tetrahedron form a regular tetrahedron, any choice of smoothing one edge among pairs of dual edges will give the same area and the same volume.

Given $M$ a Meissner polyhedron, its surface area is composed of the following elements:
\begin{itemize}[noitemsep,topsep=0pt]
	\item For every face $\tau(x)$ of $M$, having vertices $v_1,...,v_k$ consider the spherical polygon determined by $v_1,...,v_k$ on the sphere centered at the vertex $x$ opposite to $\tau(x)$. This provides a series of spherical geodesic spherical polygons.
	\item For every pair of opposite edges $(e,e')$, consider the intersection of $\partial M$ with the dihedral angles at $e, e'$ in the tetrahedron determined by the vertices of these edges. This produces two surfaces: a spherical spindle $S(e')$ (Figure \ref{fig:spindle}) and a wedge $W(e)$ (Figure \ref{fig:wedges}). The areas of these two surfaces will be established explicitly in the following in terms of the lengths $\theta(e),\theta(e')$. 
\end{itemize}

The surface area of $M$ is made of geometric surfaces which have explicit formulas for their areas. For example, it seems obvious that $S(e')$ and $W(e)$ can be computed in terms of the spherical lengths of edges $e,e'$, given by $\theta(e),\theta(e')$. Nevertheless, it is not clear for the moment how the area of the portions of the sphere present in the faces of $M$ could be computed in terms of the same parameters. This will be addressed in the next section.

\section{Surface areas of Meissner polyhedra}
\label{sec:area-computations}
\subsection{Linked rectangles on the sphere}
\label{sec:linked}

Let $M$ be a Meissner polyhedron given by \eqref{eq:def_meissner}. If $x_i$ is one of its vertices then the opposite face contains a spherical polygon $S(x_i)$ having common vertices with $M$, i.e. a region bounded by a sequence of geodesic arcs $\arc{v_iv_{i+1}}$ on a sphere of radius $1$, where $v_i, i=0,...,k-1$ are vertices of $M$. The classical Gauss-Bonnet formula allows to compute the area of a spherical polygon in terms of its angles, but it involves the turning angles of the spherical polygon which are not explicit in terms of $\theta(e),\theta(e')$. We show below how a simple observation can lead to an explicit formula for the total area of all such spherical regions in $\partial M$.

Consider the following definition of a normal vector for a bounded convex body $M$.

\begin{deff}
	Let $M$ be a bounded convex body in $\Bbb{R}^3$ and $x \in \partial M$ be a point on its boundary. A normal vector to $M$ at $x$ is a unit vector which is a normal vector for a supporting plane $\alpha$ for $M$ at $x$, pointing in the half-space determined by $\alpha$ which does not contain $M$. Denote by $N(x)$ the set of unit normal vectors to $M$ at $x$, identified with a subset of the unit sphere $\Bbb{S}^2$.
\end{deff}
We have the following elementary properties. 
\begin{prop}\label{prop:normals}
	a) The map $N: \partial M \to \Bbb{S}^2$ is surjective. 
	
	b) Suppose $M$ is strictly convex. Then $x \neq y$ implies $N(x) \cap N(y) = \emptyset$.
	
	c) If $M_0\subset M$ is a spherical surface of unit radius then $N(M_0)$ is a translation of $M_0$.
	
	d) If $M$ is a Meissner polyhedron and $x$ is a vertex then $N(x)$ contains the antipodal set in $\Bbb{S}^2$ of the spherical surface $N(S(x))$.
\end{prop}

\emph{Proof:} a) For every orientation there exists a supporting plane orthogonal to it, leaving $M$ on the opposite side, therefore the normal map is surjective.

b) If $x\neq y$ and $n \in N(x)\cap N(y)$ then parallel supporting planes exist at $x$ and $y$, showing that these planes should coincide, contradicting the strict convexity.

c) This property is obvious by definition.

d) If $x$ is a vertex and $S(x)$ is the spherical polygon contained in the face opposite to $x$ then any supporting plane $\alpha$ at a point in $S(x)$ has a unique corresponding parallel supporting plane at $x$. The desired conclusion follows.
\hfill $\square$



We arrive at the following result:

\begin{prop}\label{prop:partition}
	Let $M$ be a Meissner polyhedron. The unit sphere $\Bbb{S}^2$ can be partitioned in the following regions:
	
	(A) Translations of the spherical polygons $S(x_i)$ contained in the faces opposite to $x_i$ of $M$ and their antipodals in $\Bbb{S}^2$ which contain normals to $x_i$.
	
	(B) A series of pairs of antipodal rectangles having spherical edge lengths $(\theta(e_i),\theta(e_i'))$, $i=1,...,m-1$, where $(e_i,e_i')$ are the pairs of dual edges in $M$.
\end{prop}

\emph{Proof:} Following results of Proposition \ref{prop:normals}, the images of $S(x_i)$ through the normal map $N$ are spherical geodesic polygons which are translations of $S(x_i)$. Every normal to $S(x_i)$ has an opposite one at $x_i$. 

Let us now consider the complentary region to $\bigcup_{i=1}^{m-1} (N(S(x_i))\cup -N(S(x_i)))$ in $\Bbb{S}^2$. Given a pair of dual edges $(e,e')$ of $M$, with $x,y$ and $x',y'$ endpoints of $e,e'$, respectively we have the following observations. 

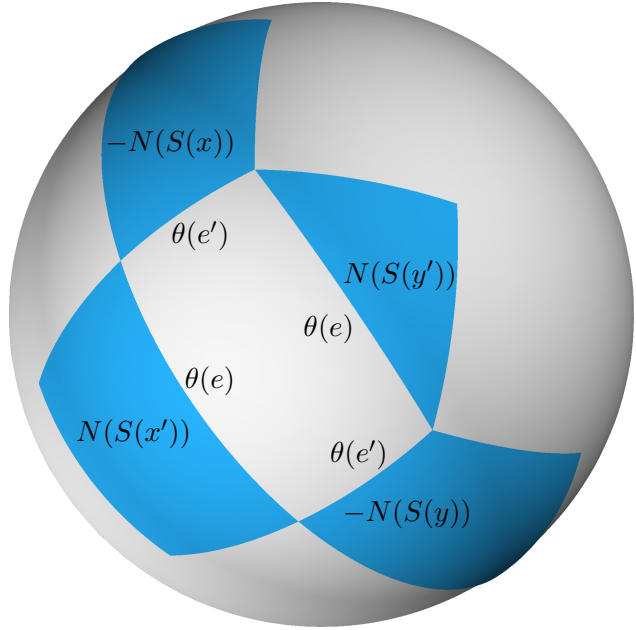
\begin{figure}
	\centering 
	\begin{picture}(200,200)
	\put(0,0){	\includegraphics[width=0.5\textwidth]{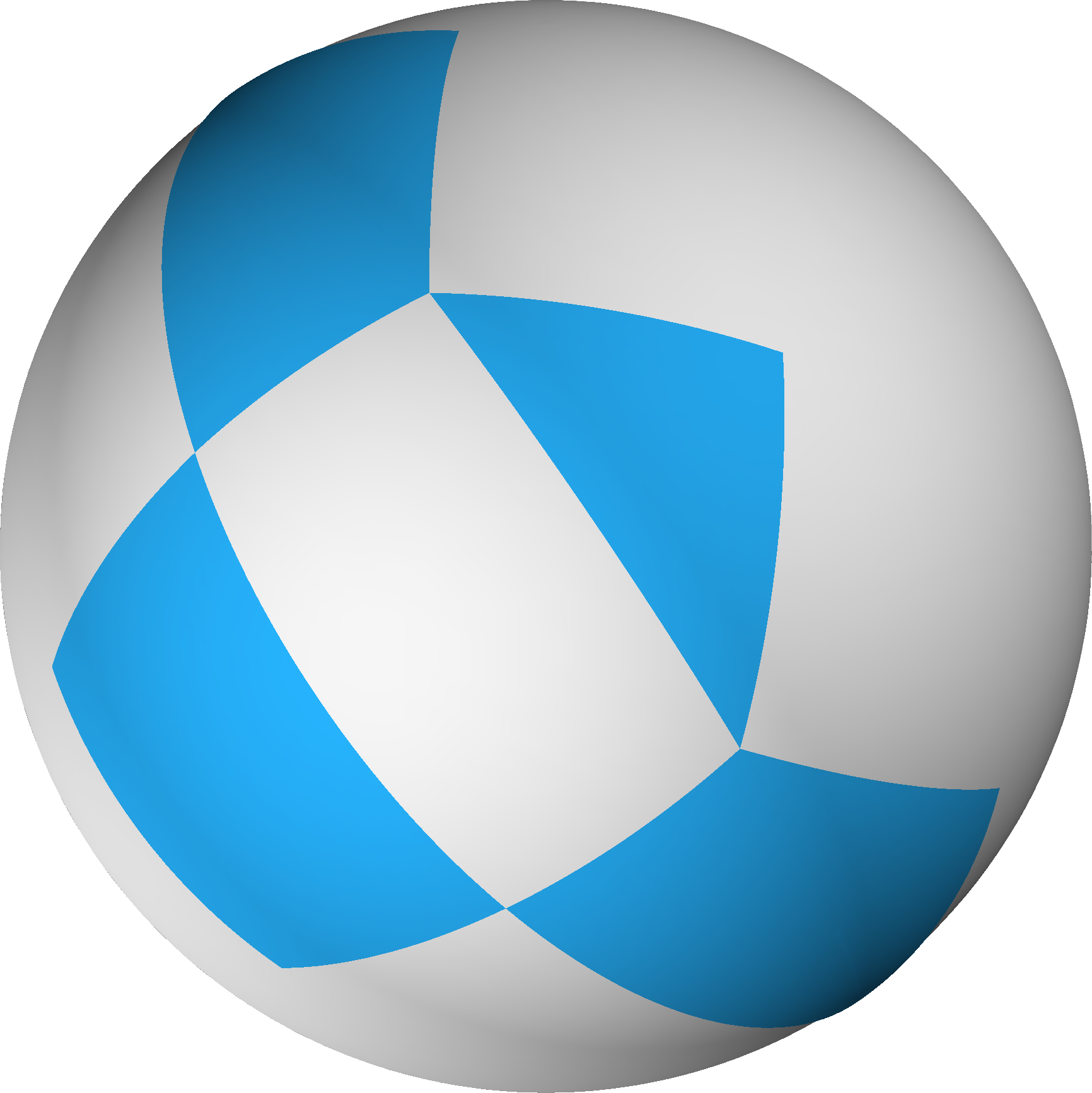}}
	\put(29,70){\Label{$N(S(x'))$}}
	\put(130,130){\Label{$N(S(y'))$}}
	\put(40,180){\Label{$-N(S(x))$}}
	\put(130,40){\Label{$-N(S(y))$}}
	\put(70,90){\Label{$\theta(e)$}}
	\put(115,110){\Label{$\theta(e)$}}
	\put(65,145){\Label{$\theta(e')$}}
	\put(125,62){\Label{$\theta(e')$}}
	\end{picture}
	\caption{Example of spherical rectangle corresponding to normals for a wedge or spindle surface around the edge $e$ with endpoints $x,y$ and dual to $e'$ with endpoints $x', y'$. The adjacent spherical polygons to $e$ are $S(x')$, $S(y')$. The spherical region bounded by $N(S(x')), N(S(y')), -N(S(x)), -N(S(y))$ is a spherical rectangle with edge lengths $\theta(e),\theta(e')$.  
	}
	\label{fig:rectangles-labels}
\end{figure}

Normals considered in part (A) for vertices $x,y$ and faces $\tau(x'),\tau(y')$ adjacent to the edge $\arc{xy}$ have the following representations in $\Bbb{S}^2$:
\begin{itemize}[topsep=0pt,noitemsep]
	\item the spherical polygons $N(S(x')), N(S(y'))$, having an edge equal to $\theta(e)$ corresponding to the image through $N$ of a geodesic from $x$ to $y$.
	\item the spherical polygons $-N(S(x)),-N(S(y))$ having an edge equal to $\theta(e')$ corresponding to the image through $N$ of a geodesic from $x'$ to $y'$.
\end{itemize}
An illustration is given in Figure \ref{fig:rectangles-labels}.
Therefore, normals to $W(e)$ or $S(e)$ not considered in part (A) are contained in a spherical rectangle having edges $\theta(e),\theta(e')$. Indeed, this region is a spherical quadrilateral with equal opposite sides and two planes of symmetry. Moreover, its vertices form an Euclidean rectangle. This is a consequence of the fact that spindle and wedge surfaces also have two planes of symmetry. This determines a \emph{spherical rectangle}, i.e. a spherical quadrilateral with equal opposite sides and equal angles.

The same configuration arises for normals around $e'$ not considered at part (A). The resulting spherical rectangle is the antipodal for the one obtained for $e$. \hfill $\square$

\begin{figure}[!htp]
	\centering
	\includegraphics[width=0.3\textwidth]{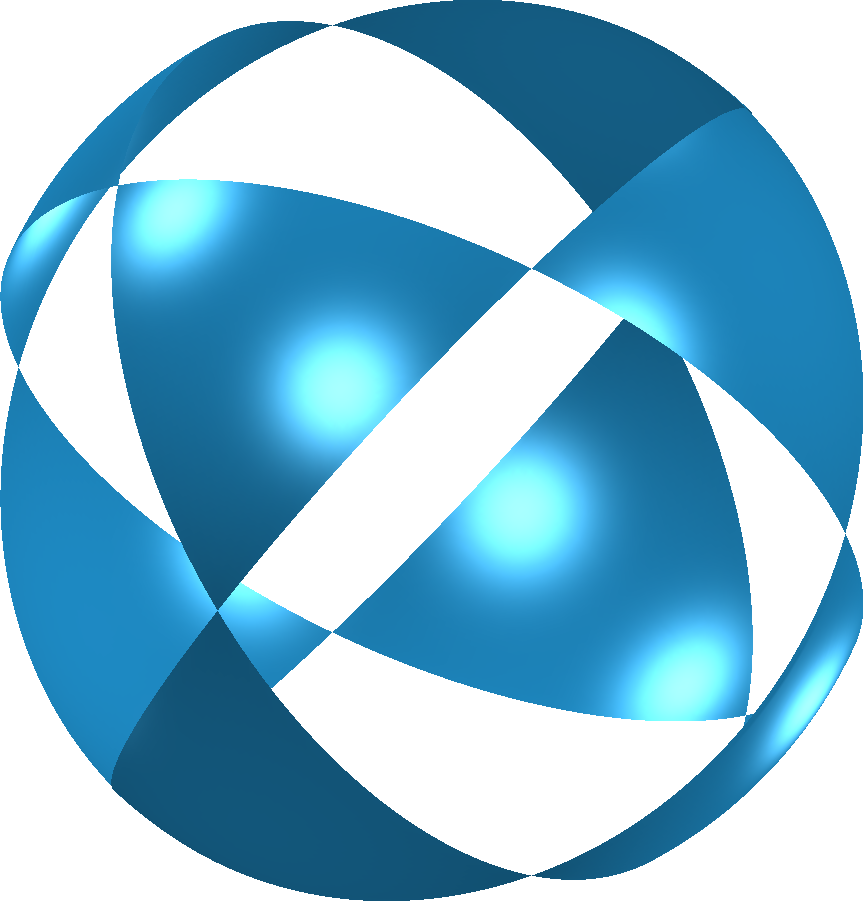}\quad
	\includegraphics[width=0.3\textwidth]{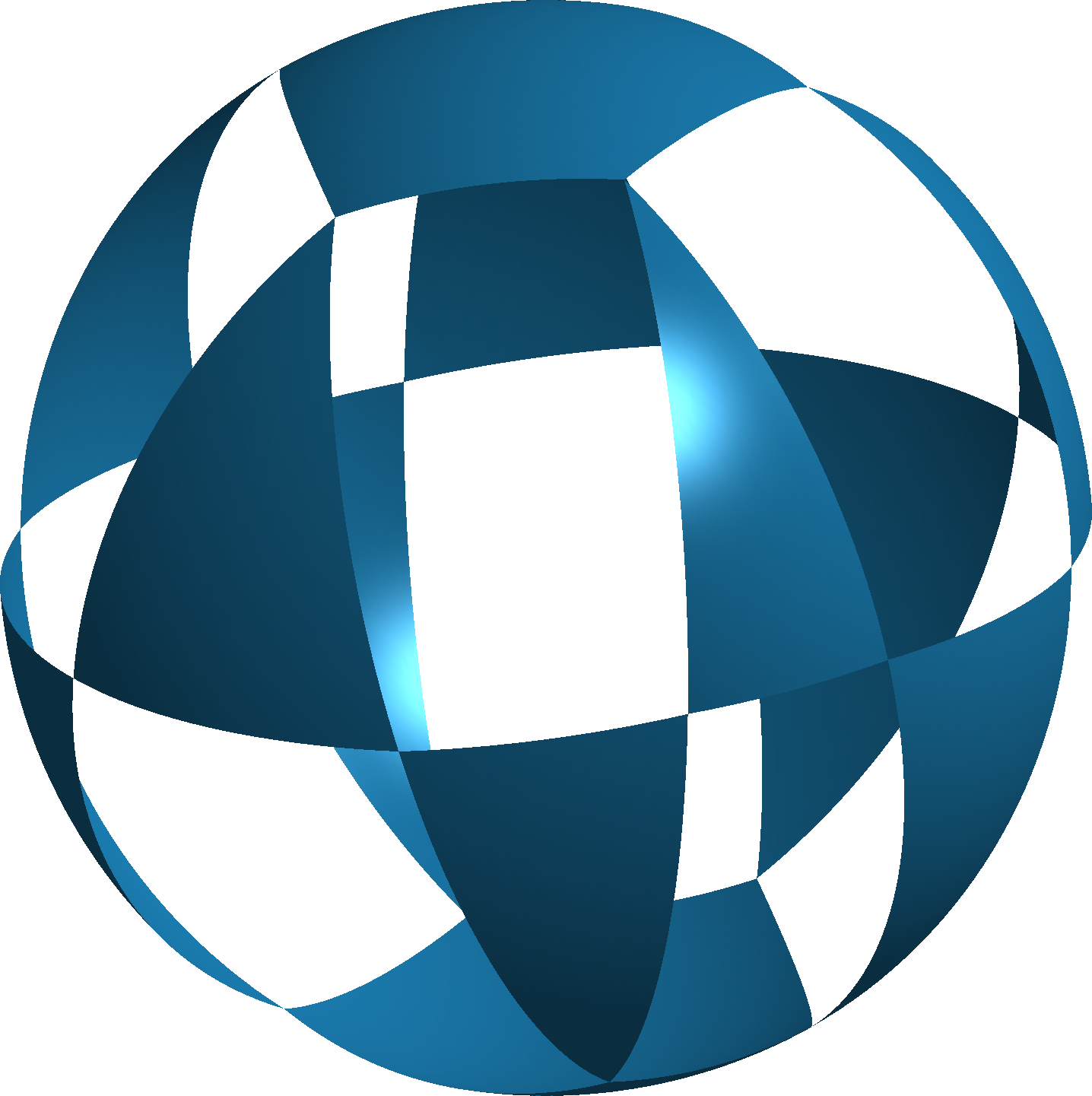}\quad
	\includegraphics[width=0.3\textwidth]{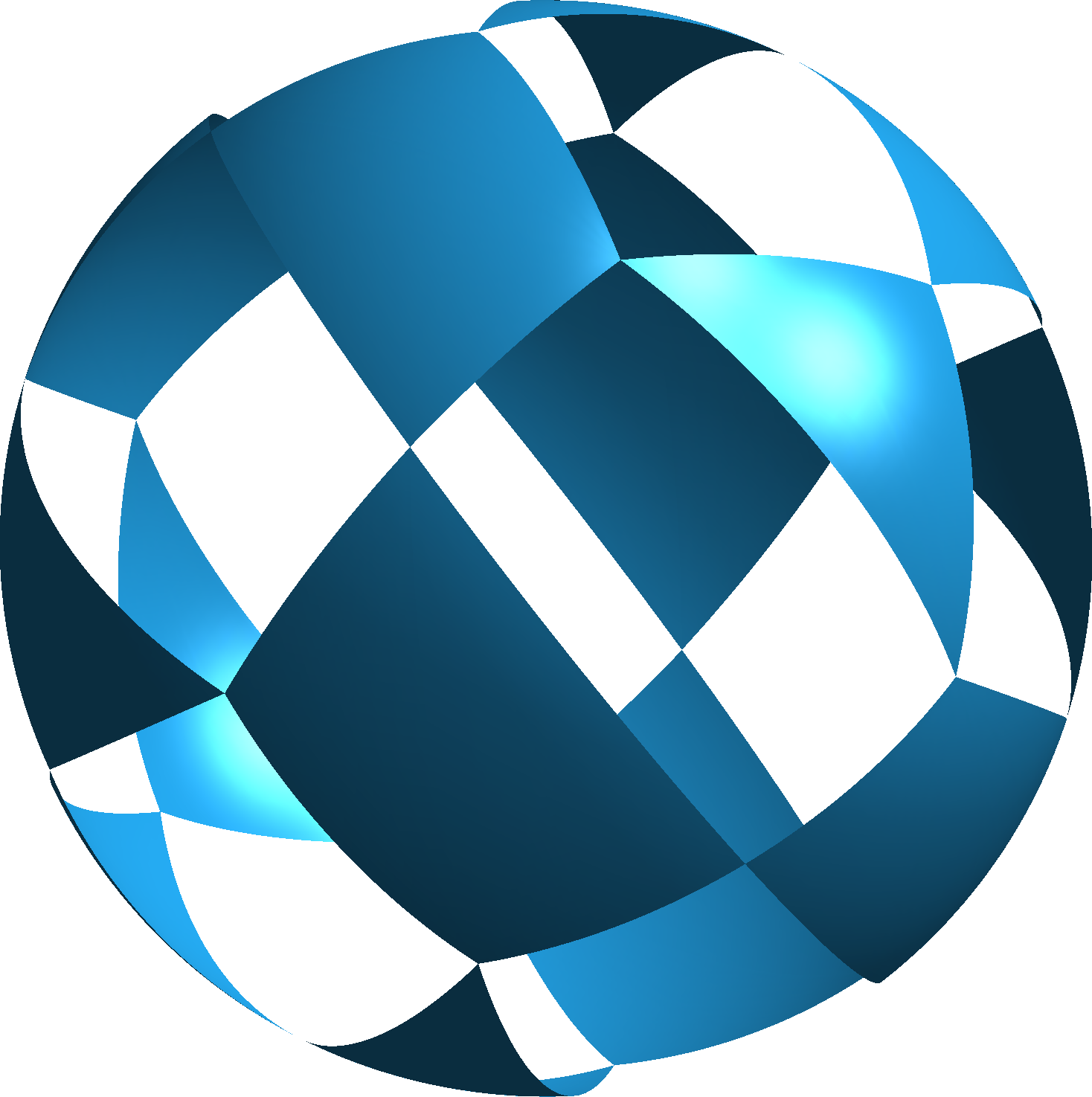}\quad
	\caption{Examples of partitions of the sphere containing spherical parts of Meissner or Reuleaux polyhedra and spherical rectangles.}
	\label{fig:rectangles-sphere}
\end{figure} 

To compute the surface area of a Meissner polyhedron $M$ a clear strategy emerges, following the previous results. The area of spherical regions in $\partial M$ is the complementary of a series of spherical rectangles with lengths $(\theta(e_i),\theta(e_i'))_{i=1}^{m-1}$. The remaining parts are wedge and spindle surfaces whose area will be computed. We start by computing the area of a spherical rectangle.

\begin{lemma}\label{lem:rectangles}
A spherical rectangle is a spherical quadrilateral with equal angles and equal opposite sides. The area of a spherical rectangle with geodesic side lengths $\theta, \theta'\in (0,\pi)$ is given by
\[ R(\theta,\theta')= 4\arcsin \left(\tan \frac{\theta}{2} \tan \frac{\theta'}{2}\right).\]
\end{lemma}

\begin{figure}
	\centering
	\includegraphics[width=0.5\textwidth]{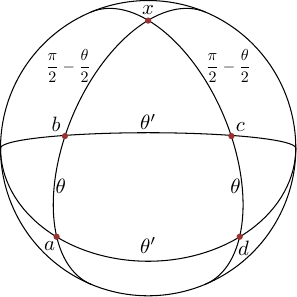}
	\caption{Computing the area of a spherical rectangle in terms of its side lengths.}
	\label{fig:rectangle-area}
\end{figure}

\emph{Proof:} Let $a,b,c,d$ be the vertices of the spherical rectangle such that $\arc{ab}, \arc{cd}$ have lengths $\theta$ and $\arc{bc},\arc{da}$ have lengths $\theta'$. To compute the area of a spherical geodesic polygon, its angles need to be computed. See Figure \ref{fig:rectangle-area} for an illustration.

Draw the great circles containing $\arc{ab},\arc{cd}$ and denote by $x$ the intersection point closer to $\arc{bc}$ and $x'$ the other point of intersection. By the symmetry of the rectangle, the arc $\arc{ab}$ is symmetric about the midpoint of the half circle $\arc{xx'}$ containing $\arc{ab}$. For simplicity, we use the same notation for an arc and for its spherical length. Thus, in the spherical triangle $xbc$ all edge lengths are known: $\arc{bc} = \theta'$, $\arc{bx}=\arc{cx} = \frac{\pi}{2}-\frac{\theta}{2}$. Denoting by $\alpha$ the angle at $b$ in the triangle $bcx$, the spherical law of cosines gives
\[ \cos \arc{cx} = \cos \arc{bc} \cos \arc{bx}+\sin \arc{bc}\sin \arc{bx} \cos \alpha.\]
Therefore, using basic trigonometric identities, we obtain
\[ \cos \alpha = \frac{\cos \left(\frac{\pi}{2}-\frac{\theta}{2}\right)-\cos \theta' \cos\cos \left(\frac{\pi}{2}-\frac{\theta}{2}\right)}{\sin \theta' \sin \left(\frac{\pi}{2}-\frac{\theta}{2}\right)}=\frac{\sin \frac{\theta}{2}(1-\cos \theta')}{\sin \theta'\cos \frac{\theta}{2}}=\tan \frac{\theta}{2}\tan \frac{\theta'}{2}. \]
The Gauss-Bonnet formula combined with the result above shows that the area of the rectangle is
\[ R(\theta,\theta') = 2\pi-4\alpha=4\left(\frac{\pi}{2}-\alpha\right)=4\arcsin \left(\tan \frac{\theta}{2}\tan \frac{\theta'}{2}\right).\]
\hfill $\square$

The area of $\Bbb{S}^2$ is $4\pi$. Therefore, keeping in mind that every spherical part counted in (A) in Proposition \ref{prop:partition} and every rectangle appears twice on the sphere in Proposition \ref{prop:partition}, it follows that the total area of the spherical polygons contained in $\partial M$ is 
\[ 2\pi-\sum_{i=1}^{m-1} R(\theta(e_i),\theta(e_i')).\]

Thus, we arrive at an initial formula for the area of a Meissner polyhedron, which is given by:
\begin{equation}\label{eq:area-Meissner-v1}|\partial M| = 2\pi-\sum_{i=1}^{m-1} R(e_i,e_i')+\sum_{i=1}^{m-1}( |W(e_i)|+|S(e_i')|),
\end{equation}
where $|W(e_i)|,|S(e_i')|$ represent the areas of \emph{wedge} and \emph{spindle} regions, respectively. 

Thus, the Blascke-Lebesgue problem in dimension three amounts to solving
\begin{equation}\label{eq:BL-reformulation}
 \max \sum_{i=1}^{m-1} \left(R(e_i,e_i')-|W(e_i)|-|S(e_i')|.\right).
 \end{equation}
In the following section $|W(e_i)|$ and $|S(e_i')|$ are computed explicitly in terms of $\theta(e_i),\theta(e_i')$.

\subsection{Computation of areas of wedge and spindle surfaces}

Computations regarding spindle and wedge surfaces are also presented \cite{hynd-vol-per}, \cite{meissner_hynd}. The computations made below sometimes use different arguments, therefore we present the computations in full detail, for the sake of completeness.

In the following $(e,e')$ denotes a generic pair of dual edges and $\theta(e),\theta(e')$ denote their spherical lengths. Recalling that $\phi(e)$ is the dihedral angle  at edge $e$ of the tetrahedron determined by $e$ and $e'$, we have the relation 
\begin{equation}\label{eq:dihedral} \sin \frac{\phi(e)}{2} = \frac{\sin \frac{\theta(e')}{2}}{\cos \frac{\theta(e)}{2}}.
\end{equation}
Since $\theta(e) \in [0,\pi/3]$ for every edge, we find that $\sin \frac{\phi(e)}{2} \leq \frac{\sin \frac{\pi}{6}}{\cos \frac{\pi}{6}} = \frac{\sqrt{3}}{3}$. This gives $\cos \phi(e) = 1-2\sin^2 \frac{\phi(e)}{2} \geq \frac{1}{3}$ and as a consequence all dihedral angles $\phi(e)$ verify
\begin{equation}\label{eq:bound-dihedral}
\phi(e) \in [0,\arccos\frac{1}{3}].
\end{equation}

Moreover, denoting $d(e,e')$ the distance between the midpoints of the segments determined by vertices of the dual pair $(e,e')$ we have
\begin{equation}\cos \frac{\phi(e)}{2} = \frac{d(e,e')}{\cos \frac{\theta(e)}{2}}.
\label{eq:angles-vs-d}
\end{equation}
Since $d(e,e')$ depends only on $\phi(e)$, $\theta(e)$ via \eqref{eq:angles-vs-d} it follows that
\begin{equation} \cos \frac{\phi(e)}{2}\cos \frac{\theta(e)}{2}=\cos \frac{\phi(e')}{2}\cos \frac{\theta(e')}{2}.
\label{eq:angles-vs-d2}
\end{equation}

The wedge $W(e)$ (see Definition \ref{def:wedge}) is the union of two spherical regions contained between two circles. Gauss-Bonnet formula can be used to compute its surface area, provided the radii of the two circles and the angle made by the two circles are known.

{\bf Angle made by two circles on the sphere.} A circle on the unit sphere is determined by a point $C$ on the sphere and a spherical distance $\theta \in [0,\pi/2]$ and consists of all points $X$ on the sphere at spherical (or angular) distance $\theta$ to $P$. A spherical circle is, of course, an Euclidean circle of radius $\sin \theta$.

\begin{lemma}	\label{lemma:angle-circle}
	Consider intersecting circles on the unit sphere having centers $C_1,C_2$ and spherical radii $\theta_1, \theta_2 \in [0,\pi/2]$. Suppose that their axes of symmetry make an angle equal to $\phi$. Then, denoting by $\alpha$ the angle made by the two circles at one intersection point, we have
	\[ \cos \phi = \cos \theta_1\cos \theta_2+\sin \theta_1\sin \theta_2 \cos \alpha.\]
\end{lemma}

\begin{figure}
	\centering 
	\includegraphics[width=0.4\textwidth]{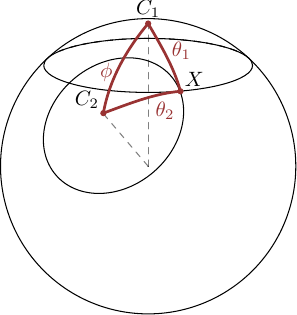}
	\caption{Configuration for computing the angle between two circles on the sphere. The two circles have centers $C_i$ and spherical radii $\theta_i$, $i=1,2$ and meet at point $X$.}
	\label{fig:circle-angle}
\end{figure}

\emph{Proof:} An illustration of the configuration is given in Figure \ref{fig:circle-angle}. Assume the two circles intersect: a necessary and sufficient condition is that $\phi$, $\theta_1$, $\theta_2$ are the sides of a spherical triangle, i.e., they verify the usual triangular inequalities. Let $X$ be a point of intersection of the two circles. Then  the arcs $\arc{C_1X}, \arc{C_2X}$ are orthogonal to the tangent vectors to the two circles at $X$. The angle of the two tangent vectors is the angle $X$ in the spherical triangle $XC_1C_2$. The sides of this spherical triangle are known, starting from the hypothesis: $C_1X=\theta_1, C_2X=\theta_2, C_1C_2 = \phi$. The spherical law of cosines coincides with the desired formula, since the spherical angle $\alpha = \angle C_1XC_2$ is opposite to the side $C_1C_2$. \hfill $\square$

{\bf Area of a wedge surface.} A wedge surface is twice the area of the intersection of a circle of a given radius smaller than $1$ and a great circle of the same sphere of radius $1$. It can be explicited using Lemma \ref{lemma:angle-circle}.

For a pair of dual edges $(e,e')$ consider the endpoints $x,y$ of $e$ and $x',y'$ of $e'$. Denoting $\arc{xy}$ the small arc between $x$ and $y$ in $\partial B(x')\cap \partial B(y')$ let us find the angle made by $\arc{xy}$ with the geodesic linking $x$ and $y$ on $\partial B(x')$ (or $\partial B(y')$. Let us identify all elements used in Lemma \ref{lemma:angle-circle}.

\begin{itemize}[noitemsep,topsep=0pt]
	\item The arc $\arc{xy}$ lies on the circle $\partial B(x')\cap \partial B(y')$ which has spherical radius $\pi/2-\theta(e')$.
	\item The geodesic lies on a big circle of spherical radius $\theta_2=\pi/2$.
	\item The angle $\phi$ made by the symmetry axes of the two circles is equal to the dihedral angle of the planes of the two circles. This is half of the dihedral angle at edge $e=xy$ in the tetrahedron $x,y,x',y'$. Indeed, the arc $\arc\{xy\}$ is contained in the bisector plane of the tetrahedron at edge $xy$, while the geodesic is contained in one of the planes $x'xy$ or $y'xy$. Thus $\phi = \phi(e)/2$.
\end{itemize}

Plugging this information into Lemma \ref{lemma:angle-circle} we find that the corresponding angle $\alpha(e)$ between the arc $\arc{xy}$ and a geodesic arc at $e$ verifies
\begin{equation}\label{eq:alpha(e)}
 \cos \alpha(e) \cos \frac{\theta(e')}{2} = \cos \frac{\phi(e)}{2}.
\end{equation}
In view of \eqref{eq:angles-vs-d}, \eqref{eq:angles-vs-d2}, notice that this angle is the same for the two pairs of opposite edges, i.e. $\alpha(e)=\alpha(e')$. It should be noted that \eqref{eq:alpha(e)} has solutions since \eqref{eq:dihedral} implies $\phi(e) \geq \theta(e')$, therefore $\cos \frac{\phi(e)}{2} \leq \cos \frac{\theta(e')}{2}$.

For a spherical region $p_1,...,p_k$ bounded by arcs of circle  on the unit sphere, having radii $r_1,...,r_k \in (0,1]$, denoting $\gamma_1,...,\gamma_k$ the turning (exterior angles) at the respective vertices, the Gauss-Bonnet formula states that its area is given by
\[ A = 2\pi - \sum_{i=1}^k \gamma_i - \sum_{i=1}^k \int_{p_i} \frac{\sqrt{1-r_i^2}}{r_i}.\]
The last term corresponds to the integral of the geodesic curvature, equal to $\frac{\sqrt{1-r^2}}{r}$ for a circle of radius $r$. 

One half of the wedge $W(e)$ is determined by two arcs of circles, meeting at an angle $\alpha(e)$ given by \eqref{eq:alpha(e)}. The geodesic arc has zero geodesic curvature. The arc $\arc{xy}$ belongs to a circle of radius $\cos \frac{\theta(e')}{2}$, has geodesic curvature $\tan \frac{\theta(e')}{2}$ and has an angular measure equal to $\phi(e')$ (see Figure \ref{fig:dual-edges}). We arrive, thus at the following result. 

\begin{prop}\label{prop:wedge-area}
	The area of a wedge surface $W(e)$ at the edge $e$ in the dual pair $(e,e')$ is given by
	\[ |W(e)| =  4\arccos\left(\frac{\cos \frac{\phi(e)}{2}}{\cos \frac{\theta(e')}{2}}\right) -2\sin \frac{\theta(e')}{2}\phi(e').\]
\end{prop}

\emph{Proof:} The proof is immediate, noting that the turning angles at the two vertices of $W_e$ are equal to $\pi-\alpha(e)$, where $\alpha(e)$ verifies \eqref{eq:alpha(e)}. Furthermore, the integral of the geodesic curvature on $\arc{xy}$ gives the second term. The result is multiplied by two, since $W(e)$ contains two congruent spherical regions. \hfill $\square$


The area of a spindle surface at edge $e$ (see Definition \ref{def:spindle}) is computed in \cite{meissner_hynd}. Nevertheless, since its computation is immediate via integrals on surfaces of revolution, we give a brief proof below.

\begin{prop}\label{prop:spindle-area}
 The area of a spindle surface $S(e)$ at the edge $e$ in the dual pair $(e,e')$ is given by
\[ |S(e)|=2\phi(e)\left(-\cos \frac{\theta(e)}{2}\frac{\theta(e)}{2}+\sin \frac{\theta(e)}{2}\right).\]
\end{prop}

\emph{Proof:} We compute the full surface area of a complete spindle and recover the result for a spindle of angle $\phi(e)$ as a byproduct. Recall that if $f : [a,b]\to \Bbb{R}_+$ is a $C^1$ function then the area of the surface of revolution determined by the graph of $f$ around the $x$-axis is simply $A = 2\pi \int_a^b f(x)\sqrt{1+(f'(x))^2)}dx$.

We apply this result for the function $f: [-\sin \frac{\theta(e)}{2},\sin \frac{\theta(e)}{2}]$ which verifies 
\[ x^2+\left(f(x)+\cos \frac{\theta(e)}{2}\right)^2 = 1.\]
We obtain
\[ f(x) = \sqrt{1-x^2}-\cos \frac{\theta(e)}{2},\ f'(x) = \frac{-x}{\sqrt{1-x^2}}.\]
The result follows immediately from the integral formula. \hfill $\square$

Notice that if in the definition of the Meissner polyhedron $M$ we consider balls with centers on $e$, then close to edges in the dual pair $(e,e')$ the boundary $\partial M$ is made of the wedge $W(e)$ and the spindle surface $S(e')$. It is already apparent from Propositions \ref{prop:wedge-area}, \ref{prop:spindle-area} that when computing $|W(e)|+|S(e')|$ some terms cancel. More precisely, 
\[ |W(e)|+|S(e')|=  4\alpha(e) -\phi(e')\theta(e')\cos \frac{\theta(e')}{2}.\]

Moreover, using the classical identity $\arcsin x=\arccos \sqrt{1-x^2}$ for $x=\tan \frac{\theta(e)}{2}\tan \frac{ \theta(e')}{2} \in [0,1]$, we find, remarkably, that 
\[\arcsin\left( \tan \frac{\theta(e)}{2} \tan \frac{\theta(e')}{2}\right)=\arccos\left(\frac{\sqrt{\cos^2 \frac{\theta(e)}{2}-\sin^2 \frac{\theta(e')}{2}}}{\cos \frac{\theta(e)}{2}\cos \frac{\theta(e')}{2}}\right)=\arccos\Big(\frac{\cos \frac{\phi(e)}{2})}{\cos \frac{\theta(e')}{2}}\Big).\]
Therefore, $4\alpha(e)$ coincides with the area of the spherical rectangle with edges $(\theta(e),\theta(e'))$, which further simplifies the expression \eqref{eq:area-Meissner-v1}. 

\subsection{Area of Meissner polyhedra}

Gathering all the results from previous sections we obtain the following.
\begin{thm}\label{thm:surface-Meissner}
	The Meissner polyhedron $M = B(X\cup e_1\cup...\cup e_{m-1})$ given by \eqref{eq:def_meissner} has the surface area
	\begin{equation}\label{eq:area-Meissner}
	|\partial M | = 2\pi - \sum_{i=1}^{m-1} f(\theta(e_i),\theta(e_i'))
	\end{equation}
	with 
	\begin{equation}\label{eq:def_f}
	f(x,y) = y\cos \frac{y}{2} \arcsin \left( \frac{\sin \frac{x}{2}}{\cos \frac{y}{2}}\right) 
	\end{equation}
\end{thm}

The maximal value of $f(x,y)$ given in \eqref{eq:def_f} for $\{(x,y) \in [0,\pi/3]^2: x \leq y\}$ is attained at $(\pi/3,\pi/3)$ and $f$ is increasing in both $x,y$ and is convex in the $x$ and $y$ directions. Results regarding the function $f$ are gathered in Appendix \ref{appA}.  

\begin{rem}\label{rem:geom-alternative}
	We have the following geometric interpretation for $f(\theta(e_i),\theta(e_i'))$. Observe that $\cos \frac{\theta(e')}{2} \phi(e')$ is the length $\ell(e)$ of the circle arc between $x, y$ in $B(x')\cap B(y')$. Therefore $f(\theta(e_i),\theta(e_i')) = \ell(e)\theta(e')$.
\end{rem}

We obtain the following direct consequence of Theorem \ref{thm:surface-Meissner}.

\begin{cor}\label{cor:consequences}
	(a) The area of a Meissner tetrahedron is equal to \[2\pi-3f(\pi/3,\pi/3)= 2\pi-\frac{\sqrt{3}}{2}\arccos \frac{1}{3} \pi.\]
	The volume of a Meissner polyhedron is obtained using the Blascke formula \eqref{eq:blaschke}. 
	
	(b) Among the $2^{m-1}$ possible choices of Meissner polyhedra for a given extremal set of diameter one consisting of $m$ points, the one with minimal area verifies
	\[ \theta(e_i)\leq \theta(e_i'),\  i = 1,...,m-1.\]
	In other words, the longest edge among each pair of dual edges should be \emph{smoothed}, i.e. replaced by a spindle surface.
	
	(c) The three dimensional Blaschke-Lebesgue theorem is equivalent to solving
	\[\max \sum_{i=1}^{m-1} f(\theta(e_i),\theta(e_i')),\]
	assuming that in all dual pairs we have $\theta(e_i) \leq \theta(e_i')$.
\end{cor}

\emph{Proof:} (a) The Meissner tetrahedron corresponds to $m=3$ and all edges have spherical length $\pi/3$. The result follows after evaluating $f(\pi/3,\pi/3)$.

(b) We observe that on $[0,\pi/3]^2$ we always have
\begin{equation}\label{eq:ineq-xy} x \leq y \Longrightarrow f(x,y)\geq f(y,x).
\end{equation}
A proof is given in the Appendix \ref{appA}. Therefore the conclusion follows. 
 
(c) A simple consequence of the area formula \eqref{eq:area-Meissner}.
\hfill $\square$ 

\begin{rem}
	The strategy employed here can also be used to compute the area of Reuleaux polyhedra $R$, like in \cite{hynd-vol-per}. Compared to Meissner polyhedra, for any pair of dual edges $(e,e')$ the corresponding wedge surfaces $W(e),W(e')$ are present in the boundary of a Reuleaux polyhedron. Therefore, its surface area is simply
	\begin{multline}\label{eq:area-Reuleaux}|\partial R| = 2\pi-\sum_{i=1}^{m-1} R(e_i,e_i')+\sum_{i=1}^{m-1}( |W(e_i)|+|W(e_i')|)\\
	 =2\pi+\sum_{i=1}^{m-1}(4\alpha(e_i)-2\sin \frac{\theta(e_i)}{2}\phi(e_i)-2\sin \frac{\theta(e_i')}{2}\phi(e_i')),
	\end{multline}
	The resulting expression is slightly more complex than \eqref{eq:area-Meissner} since fewer terms simplify. On the other hand the expression is symmetric in $\theta(e_i), \theta(e_i')$ as expected, recalling that $\alpha(e_i)=\alpha(e_i')$.
\end{rem}

\section{Meissner pyramids}
\label{sec:pyramids}

Given an extremal three dimensional finite set of diameter $1$ having $m$ points it is possible to attach a graph structure to it: two points $x,y$ are connected through an edge if and only if $|x-y|=1$. Such a graph will be called \emph{diameter-graph} in the following. These graphs are studied in detail in \cite{meissner_graphs}, \cite{polyhedra_graphs} and \cite{self-dual-graphs}. In \cite{meissner_graphs} all possible diameter graphs for $m\leq 14$ are investigated, giving rise to many examples of Meissner polyhedra. 

\begin{figure}
	\centering 
	\includegraphics[height=0.3\textwidth]{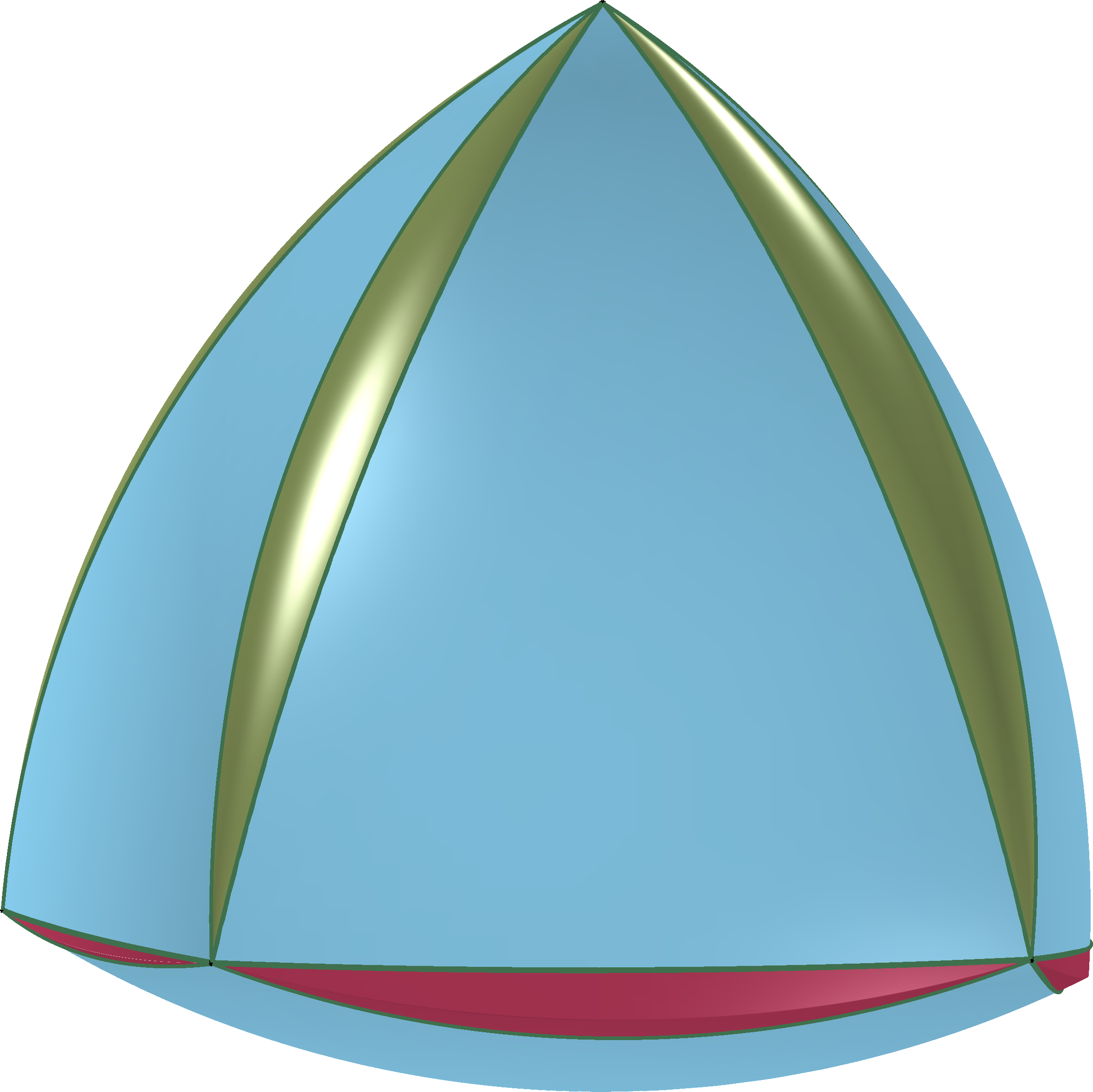}
	\includegraphics[height=0.3\textwidth]{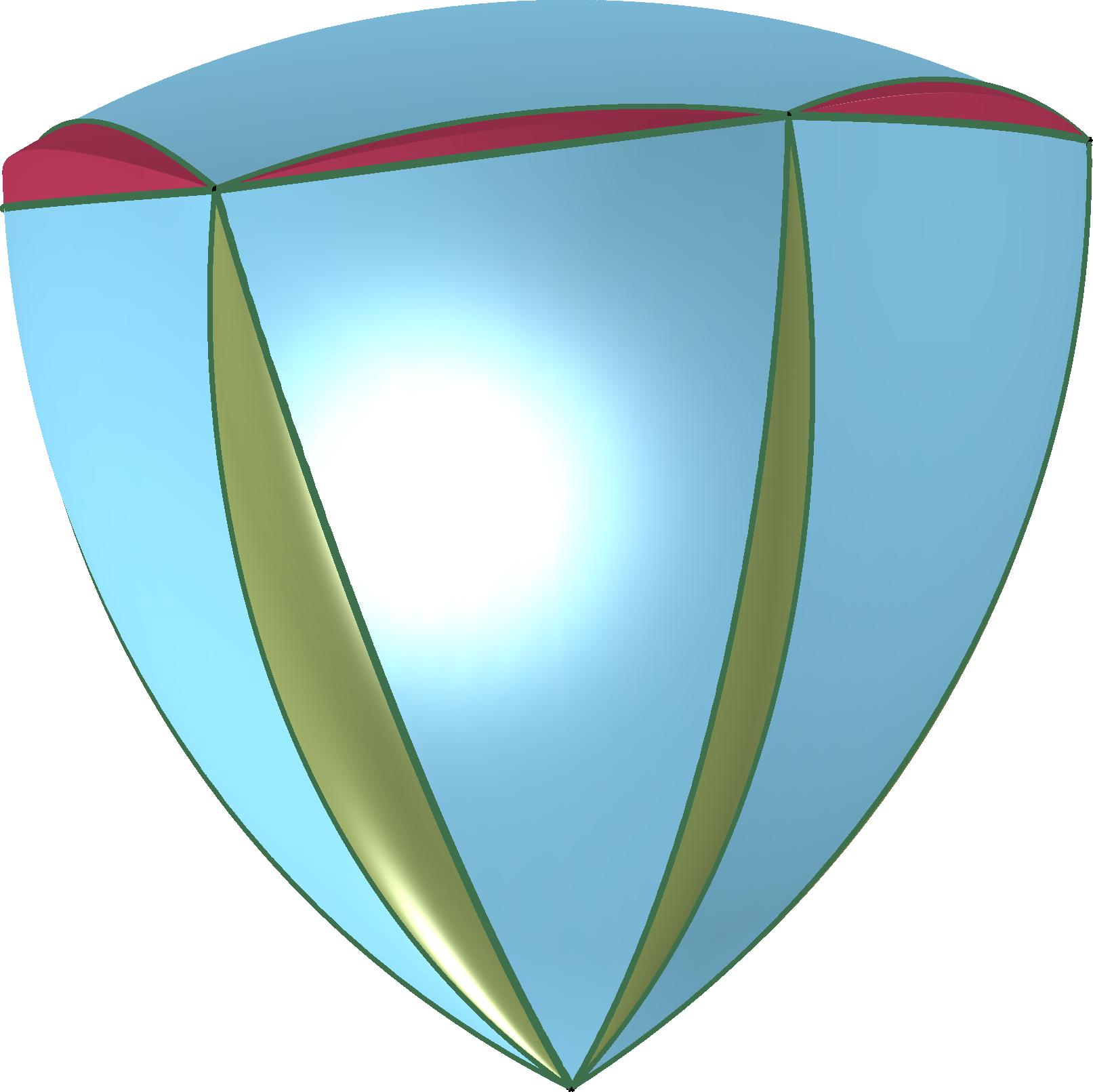}
	\includegraphics[height=0.3\textwidth]{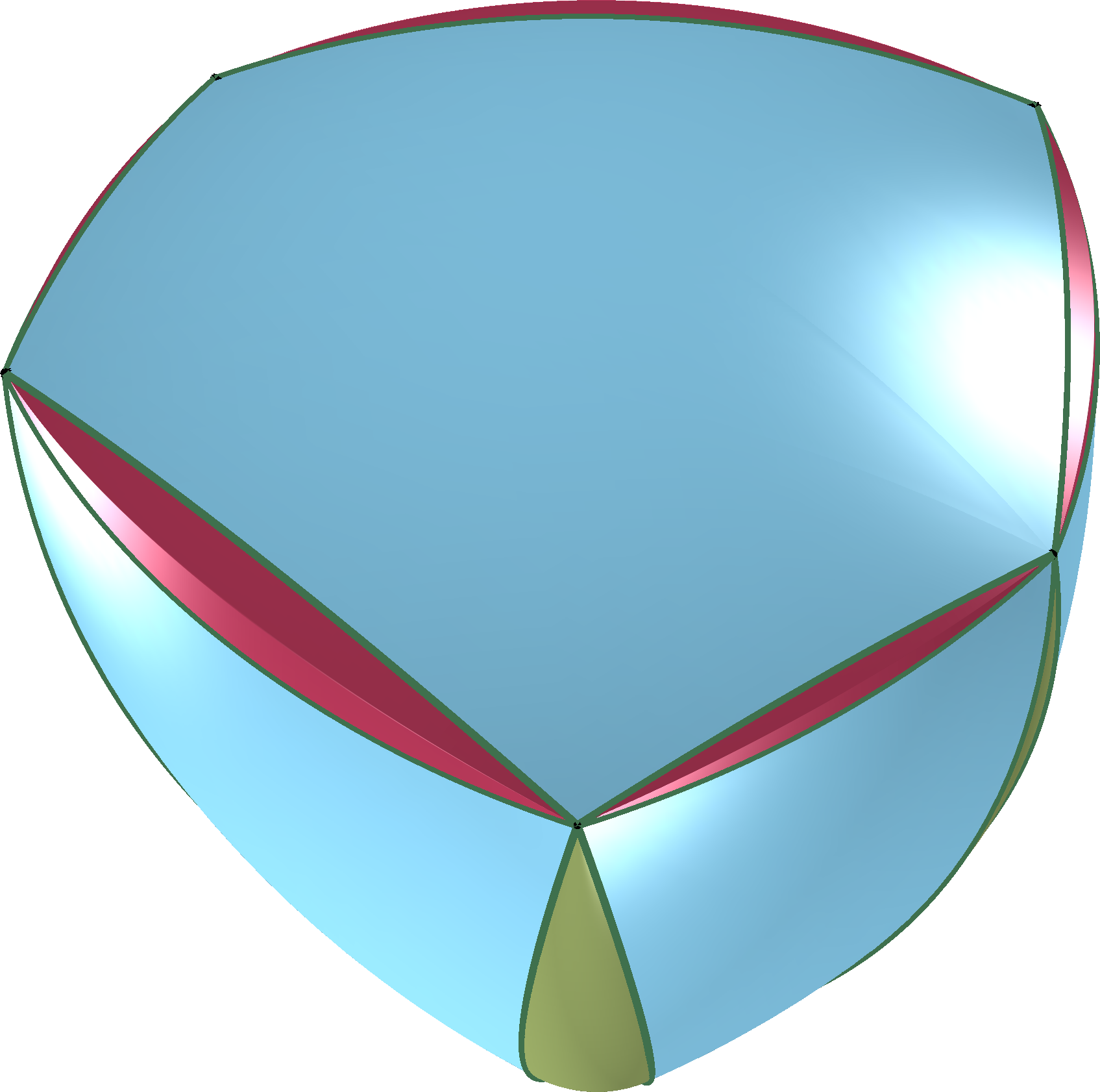}
	\caption{Different views of a Reuleaux pyramid. The edges connecting the central node $a$ to the vertices of the opposite face are smoothed. The rightmost picture shows the face opposite to the central node $a$ and is a spherical Reuleaux polygon of constant width equal to $\pi/3$ (in spherical distance).}
	\label{fig:pyramids}
\end{figure}

In \cite{polyhedra_graphs} the authors show that there is a class of graphs which is particular, in the sense that every edge is on a triangular face. Such a graph is called a \emph{wheel graph} and has one central node $a$ connected to nodes $b_1,...,b_{m-1}$ which form a cycle (see Figure \ref{fig:pyramids}). Geometrically, Meissner polyhedra having wheel diameter graphs resemble a pyramid. Such polyhedra will be called \emph{Meissner pyramids} in the following. They consist of a vertex $a$ and a face $b_1,...,b_{m-1}$ opposite to $a$. We are interested in Meissner pyramids of minimal volume, therefore the edges $ab_i$ should be smoothed. This implies that the face $b_1,...,b_{m-1}$, where wedge surfaces are considered at each one of the edges, resembles a planar Reuleaux polygon, each vertex being at Euclidean unit distance from two opposite vertices. This is, in fact, a spherical Reuleaux polygon with spherical width $\pi/3$.

In this section it is shown that the Meissner tetrahedron minimizes the surface area and volume among Meissner pyramids. 

\begin{thm}\label{thm:optim-pyramids}
	The minimal surface area of a Meissner pyramid is attained for Meissner tetrahedra. The same result holds for the volume in view of the Blaschke formula \eqref{eq:blaschke}.
\end{thm}

\emph{Proof:} Consider $M = ab_1...b_{m-1}$ a Meissner pyramid. Every pair of dual edges contains one edge of the form $b_ib_{i+1}$ and one of the form $ab_j$. Since $|a-b_j|=1$, the associated spherical length is $\theta(ab_j) = \pi/3$. According to Corollary \ref{cor:consequences}, the minimal area is attained when spindles surfaces are put along the longest edge of every one of the dual pairs. For Meissner pyramids, all edges of the form $ab_i$ have spherical length $\pi/3$ and are, therefore larger than their dual edges (of the form $b_jb_{j+1}$, with the convention $b_m=b_1$). Denoting $e_j = b_jb_{j+1}$, $j=1,...,m-1$ the least area of a Meissner pyramid having the given vertices is equal to
\[ 2\pi-\sum_{i=1}^{m-1}f(\theta(e_j),\pi/3).\]
Using the explicit expression of $f(x,y)$ in \eqref{eq:def_f} and the Remark \ref{rem:geom-alternative} regarding the alternative expression we arrive at the equivalent problem
\begin{equation}\label{eq:pyramid-max} \max \sum_{i=1}^{m-1} \frac{\pi}{3}\cos \frac{\pi}{3} \phi(ab_i) \Leftrightarrow \frac{\pi}{3}\max \sum_{i=1}^{m-1} \ell(b_ib_{i+1}), 
\end{equation}
where it should be recalled that $\phi(ab_i)$ is the dihedral angle at $ab_i$ in the tetrahedron $ab_ib_jb_{j+1}$, where $b_{j}b_{j+1}$ is the dual of $ab_i$. The faces $ab_ib_j$, $ab_ib_{j+1}$ of this tetrahedron are equilateral triangles. Also, according to Remark \ref{rem:geom-alternative}, $\ell(b_ib_{i+1})$ is the length of the arc $\arc{b_ib_{i+1}}$ on the intersection $\partial B(a)\cap \partial B(b_i)$. It is immediate to see that up to a multiplicative function, the objective function to be maximized is equal to the perimeter of the spherical region determined by the face $b_1...b_{m-1}$, including the parts coming from the wedge surfaces, on the sphere centered in $a$ with radius $1$. Moreover, all boundary parts of this region are circle arcs of radii $\frac {\sqrt{3}}{2}$, since they are on the intersection of $\partial B(a)$ and $\partial B(b_i)$, $i=1,...,m-1$.

Let us make the analogy with planar Reuleaux polygons more precise. The face opposite to $a$ in the Meissner pyramid whose spindle surfaces are all on edges containing $a$ is the intersection on $\partial B(a)$ of the spherical circles centered in $b_i$, $i=1,...,m-1$ having spherical radius $\pi/3$. This is a spherical region of constant width in $\partial B(a)$. Blaschke claimed in \cite{Blaschke1915} that the analogue of the Blaschke-Lebesgue theorem holds also on the sphere: among all shapes on a sphere of radius $1$ of constant width $w\leq \pi/2$ (in the sense of the spherical distance), the spherical Reuleaux triangle of width $w$ has the smallest area. This was proved by Leichtweiss in \cite{Leichtweiss} and strengthened to the case of spherical disk polygons in \cite{sph_disk_poly}. For spherical curves of constant width $w$ the length $L$ and the area $A$ are linked by the formula
\[ L = (2\pi-F)\tan \frac{w}{2}.\]
This formula is attributed to Blaschke \cite{Blaschke1915} and can also be found in \cite{santalo_sphere}. In particular, at fixed constant width, the spherical Reuleaux triangle maximizes the perimeter. 

In conclusion, the solution to problem \eqref{eq:pyramid-max} corresponds to a spherical Reuleaux triangle and the Meissner pyramid maximizing \eqref{eq:pyramid-max} must be a Meissner tetrahedron. \hfill $\square$

\section{The general case: discrete problems}
\label{sec:existence}
Given $m \geq 4$ denote by $\mathcal M_m$ the class of Meissner polyhedra with at most $m$ vertices. This set is non-void if $m\geq 4$, as it contains the Meissner tetrahedron. Like in the two dimensional case discussed in \cite{kupitz-martini} we may formulate a series of finite dimensional problems related to the Meissner conjecture. 

\begin{prop}\label{prop:existence}
	There exists a Meissner polyhedron $M\in \mathcal M_m$ which minimizes its surface area, i.e. the problem
	\[ \min_{M \in \mathcal M_m} |\partial M|\]
	has solutions.
\end{prop}

\emph{Proof:} The standard method in the calculus of variation is employed. Of course, $|\partial M|$ has trivial upper and lower bounds: zero and the area of the unit ball. Therefore there exists a minimizing sequence $(M_k)\subset \mathcal M_m$ such that
\[ |\partial M_k|\to \inf_{M\in \mathcal M_m} |\partial M|.\]

The classical Blaschke selection theorem \cite[Theorem 1.8.7]{schneider} shows that there exists a subsequence of $(M_k)$ which converges in the Hausdorff metric (see \cite[Chapter 2]{henrot-pierre-english}). Suppose, up to relabeling, that $M_k$ converges to $M$, which is still of constant width (this is classical, see \cite{AntunesBogosel22} for a proof). The goal is to prove that $M\in \mathcal M_m$, i.e. the class of Meissner polyhedra with a fixed number of vertices is closed in the Hausdorff metric.

Since each $M_k$ has at most $m-1$ pairs of dual edges, there is a subsequence of $(M_k)$, denoted again with $(M_k)$ which has a constant number $m_0$ of pairs of dual edges. For each of the polyhedra $M_k$ consider the pairs of dual edges
\[ (e_1^k,(e_1^k)'),...,(e_{m_0}^k,(e_{m_0}^k)').\]
Corollary \ref{cor:consequences} shows that among all $2^{m_0}$ possible choices of smoothing one of any pair of dual edges, the one verifying $\theta(e_i^k)\leq \theta((e_i^k)')$ gives the lowest surface area. If this is not the case, modify each $M_k$, eventually decreasing its surface area, recovering another minimizing sequence. At convergence continuity implies that we can assume that the choice of the edge smoothing does not change.

Up to extracting a diagonal sequence, assume that 
endpoints of $e_i^k$ and $(e_i^k)'$ converge so that the corresponding segments verify
\[e_i^k \to e_i, \ \ (e_i^k)' \to e_i'\]
for all $i=1,...,m_0$. Since endpoints of edges $e_i^k$ form extremal diameter $1$ sets, in the limiting process, these diameters are preserved. If none of the edges $e_i^k,  (e_i^k)'$ collapse in the limiting process, then the discrete set obtained is extremal, having $2m-2$ diameters. If one edge $e_i$ collapses and $(e_i^k)'$ does not then two pairs of diameters merge and one vertex is lost in the limiting process, the resulting set remaining extremal. If both edges in a pair $(e_i^k,(e_i^k)')$ collapse in the limiting process then two vertices and four diameters are lost, the remaining set still being extremal. The $m_0$ vertices of $M_n$ converge to $m_1 \geq 4$ distinct points in $\Bbb{R}^3$, since the limit set $M$ has constant width.

Since Hausdorff limits of balls are balls, and Hausdorff convergence preserves finite intersections \cite[Chapter 2]{henrot-pierre-english}, all spherical parts of $M_k$ converge to spherical parts in $M$. In particular, wedge surfaces converge to wedges.

This shows that faces of $M_k$ converge to faces of $M$ and vertices of $M_k$ converge to vertices in $M$. Since $M$ contains an extremal diameter one set with $m_1$ vertices, any edge dual to a wedge will correspond to a spindle surface \cite{meissner_hynd}. Therefore will be a Meissner polyhedron with at most $m$ vertices, i.e., $M \in \mathcal M_m$. \hfill $\square$

\begin{rem}
	An alternative proof of Proposition \ref{prop:existence} can be given using Theorem \ref{thm:surface-Meissner}. Indeed, any polyhedron in $\mathcal M_m$ can be characterized by the lengths $(\theta(e_i),\theta(e_i'))_{i=1}^{m-1}$ for the pairs of dual edges. Polyhedra with fewer than $m$ vertices can also be characterized, setting some of the lengths $\theta(e_i),\theta(e_i')$ to zero. The space of lengths of dual edges is finite dimensional, closed and blunded. According to the area formula \eqref{eq:area-Meissner} the result of Proposition \ref{prop:existence} follows since a continuous function attains its extremal values on a compact set. 
\end{rem}

\begin{rem}
	The result of Proposition \ref{prop:existence} suggests that the Meissner conjecture can be reduced to a series of discrete problems. Ideally it should be proved that a Meissner polyhedron which is not a tetrahedron cannot be optimal.
	
	A similar result can be formulated if instead of Meissner polyhedra we consider Reuleaux polyhedra. In this case, the solution is expected to change as the number of vertices increase. Indeed, it is known that Reuleaux polyhedra can approximate arbitrarily well any shape of constant width. However, since these polyhedra do not have constant width themselves, when $n \to \infty$ the solution to the problem will change and converge to a constant width body of minimal volume.
\end{rem}

\section{Concluding remarks}

In this paper, the surface area and volume of Meissner polyhedra is computed explicitly in terms of lengths of dual edges (Theorem \ref{thm:surface-Meissner}). This reduces the study of the Meissner conjecture to a series of finite dimensional problems. In particular, the minimality of the volume of the Meissner tetrahedron among Meissner pyramids is established (Theorem \ref{thm:optim-pyramids}). The answer to Conjecture \ref{conj:meissner} depends on obtaining a finer understanding of the space of lengths of dual edges $(\theta(e_i),\theta(e_i'))_{i=1}^{m-1}$. In particular, it is expected that a similar phenomenon to the one observed in dimension two occurs: no Meissner polyhedron which is not a tetrahedron can be a local minimizer for the area or the volume. 

\appendix 

\section{Analysis of a particular 2D function}
\label{appA}

The function $f(x,y) = y\cos\frac{y}{2}\arcsin \left(\frac{\sin(x/2)}{\cos(y/2)} \right)$ defined in \eqref{eq:def_f} is fundamental in understanding Meissner polyhedra with minimal surface area. In this section multiple properties are proved, with a special interest regarding the region $(x,y) \in [0,\pi/3]^2$, $x \leq y$. 

\bo{(i) $x\mapsto f(x,y), y\mapsto f(x,y)$ are increasing and convex on $[0,\pi/3]$.}

\emph{Proof:} $\sin$ is increasing on $[0,\pi/6]$ and $\arcsin$ is increasing, therefore $x \mapsto f(x,y)$ is increasing. To decide the convexity, notice that 
\begin{equation}\label{eq:dxf} \frac{\partial f}{\partial x}(x,y) = \frac{y}{2} \frac{\cos \frac{x}{2} \cos \frac{y}{2}}{\sqrt{1-\sin^2 \frac{x}{2}-\sin^2 \frac{y}{2}}}.
\end{equation}
It is straightforward to see that $x \mapsto \frac{\cos \frac{x}{2} \cos \frac{y}{2}}{\sqrt{1-\sin^2 \frac{x}{2}-\sin^2 \frac{y}{2}}}$ is increasing on $[0,\pi/2]$, therefore $x \mapsto f(x,y)$ is convex on $[0,\pi/3]$ for any fixed $y \in [0,\pi/3]$.

A straightforward computation leads to
\begin{equation}\label{eq:dyf} \frac{\partial f}{\partial y} = (\cos \frac{y}{2}-\frac{y}{2}\sin \frac{y}{2}) \arcsin \left( \frac{\sin \frac{x}{2}}{\cos \frac{y}{2}}\right)\Red{+}\frac{y}{2}\frac{\sin \frac{x}{2}\sin \frac{y}{2}}{\sqrt{1-\sin^2 \frac{x}{2}-\sin^2 \frac{y}{2}}},
\end{equation}
but it is difficult to asses the sign of this partial derivative. 

The equality \eqref{eq:dxf} shows that $f$ has the following integral representation
\begin{equation}\label{eq:integral-f}
f(x,y) = \frac{y}{2} \int_0^x  \frac{\cos \frac{t}{2} \cos \frac{y}{2}}{\sqrt{1-\sin^2 \frac{t}{2}-\sin^2 \frac{y}{2}}} dt.
\end{equation}
Therefore, since the integrand is increasing in $y$ we find that $y\mapsto f(x,y)$ is increasing. Differntiating with respect to $y$ gives
\[ \frac{\partial f}{\partial y}(x,y) = \int_0^x  \frac{\cos \frac{t}{2} \cos \frac{y}{2}}{\sqrt{1-\sin^2 \frac{t}{2}-\sin^2 \frac{y}{2}}} dt+ \frac{y}{2} \int_0^x \frac{\sin\frac{t}{2}\sin t \sin \frac{y}{2}}{4\left(1-\sin^2 \frac{t}{2}-\sin^2 \frac{y}{2}\right)^{3/2}}. \]
It is straightforward that in this formulation $y \mapsto \frac{\partial f}{\partial y}$ is increasing, and therefore $y \mapsto f(x,y)$ is convex. 

\bo{(ii) Smoothing the longest among dual edges gives a lower volume:}  
\begin{equation}\label{eq:ineqf-xy}
x\leq y \Longrightarrow f(x,y) \geq f(y,x).
\end{equation}

Differentiate \eqref{eq:dxf} with respect to $y$ to get
\begin{equation}\label{eq:dxyf}
\frac{\partial^2 f}{\partial x\partial y}(x,y) = \frac{1}{2} \frac{\cos \frac{x}{2} \cos \frac{y}{2}}{\sqrt{1-\sin^2 \frac{x}{2}-\sin^2 \frac{y}{2}}}+\frac{1}{4} \frac{y\sin \frac{y}{2} \cos \frac{x}{2} \sin^2 \frac{x}{2}}{(1-\sin^2 \frac{x}{2}-\sin^2 \frac{y}{2})^{\frac{3}{2}}}.
\end{equation}

The function 
$  x \mapsto \frac{x}{\sin x}$
is strictly increasing on $[0,\pi/3]$. Therefore, for $x\leq y$ we have
\[y\sin \frac{y}{2} \cos \frac{x}{2} \sin^2 \frac{x}{2}\geq x\sin \frac{x}{2} \cos \frac{y}{2} \sin^2 \frac{y}{2}\]
(factor $\sin^2 \frac{x}{2}\sin^2\frac{y}{2}$ on both sides) and
$$\frac{\partial^2 f}{\partial x\partial y}(x,y)\geq \frac{\partial^2 f}{\partial x\partial y}(y,x).$$
Therefore, the mapping
\[ g: x \mapsto \frac{\partial f}{\partial y}(x,y)-\frac{\partial f}{\partial x}(y,x)\]
is increasing with respect to $x$ for $x \leq y$. Since $g(0) = 0$ we have 
\[ \frac{\partial f}{\partial y}(x,y)\geq \frac{\partial f}{\partial x}(y,x)\]
whenever $0 \leq x \leq y \leq \pi/3$. Replacing $y$ by $t$ and integrating this inequality with respect to $t$ on $[x,y]$ we get $f(x,y) \geq f(y,x)$ whenever $x,y \in [0,\pi/3]$, $x \leq y$, as requested.

\bibliographystyle{abbrv}
\bibliography{./biblio}

\bigskip
\small\noindent
Beniamin \textsc{Bogosel}: Centre de Math\'ematiques Appliqu\'ees, CNRS,\\
\'Ecole polytechnique, Institut Polytechnique de Paris,\\
91120 Palaiseau, France \\
{\tt beniamin.bogosel@polytechnique.edu}\\
{\tt \nolinkurl{http://www.cmap.polytechnique.fr/~beniamin.bogosel/}}

\end{document}